\documentclass[10pt,a4paper]{amsart}
\usepackage{latexsym,graphics}

\renewcommand{\sinh}{\mathrm{sh}\,}
\renewcommand{\cosh}{\mathrm{ch}\,}

\def\LaTeX{L\kern-.36em\raise.3ex\hbox{a}\kern-.15em
    T\kern-.1667em\lower.7ex\hbox{E}\kern-.125emX}

\def\aa{{\mathcal A}}

\def\cc{{\mathcal C}}

\def\hh{{\mathcal H}}

\def\lll{{\mathcal L}}
\def\mm{{\mathcal M}}
\def\rr{{\mathcal R}}
\def\ss{{\mathcal S}}

\def\tq{\ :\ } 
\def\ffi{\varphi}
\def\eps{\varepsilon}

\def\A{{\mathbb{A}}}
\def\B{{\mathbb{B}}}
\def\C{{\mathbb{C}}}
\def\hil{{\mathbb{H}}}
\def\M{{\mathbb{M}}}
\def\N{{\mathbb{N}}}
\def\P{{\mathbb{P}}}

\def\R{{\mathbb{R}}}
\def\S{{\mathbb{S}}}

\def\X{{\mathbb{X}}}
\def\Y{{\mathbb{Y}}}
\def\Z{{\mathbb{Z}}}

\newcommand{\norm}[1]{{\left\|{#1}\right\|}}
\newcommand{\ent}[1]{{\left[{#1}\right]}}
\newcommand{\abs}[1]{{\left|{#1}\right|}}
\newcommand{\scal}[1]{{\left\langle{#1}\right\rangle}}

\newenvironment{theorem}[1][]{\noindent\sl\textbf{Theorem\ #1.}\ }{\rm}
\newenvironment{theoremnom}[2][]{\noindent\sl\textbf{Theorem #1\
(\textit{#2}).}\ }{\rm}
\newenvironment{definition}[1][]{\noindent\sl\textbf{Definition.}\
}{\rm\vskip3pt}
\newenvironment{definitionfr}[1][]{\noindent\sl\textbf{Definition.}\
}{\rm\vskip3pt}

\renewenvironment{proof}[1][]{\noindent\textit{#1}.\ \rm}{}
\newenvironment{proposition}[1][]{\noindent\sl\textbf{Proposition\ #1.}\
}{\rm\vskip3pt}

\newenvironment{lemmanum}[1][]{\noindent\sl\textbf{Lemma\ #1.}\
}{\rm\vskip3pt}
\newenvironment{lemmenum}[1][]{\noindent\sl\textbf{Lemma\ #1.}\
}{\rm\vskip3pt}
\newenvironment{lemmenom}[2][]{\noindent\sl\textbf{Lemma\ #1
(\textit{#2}).}\ }{\rm\vskip3pt}
\newenvironment{notation}[1][]{\noindent\rm\textit{Notation}\ :\
}{\rm\vskip2pt}
\newenvironment{remark}[1][]{\noindent\rm\textit{Remark}\ :\
}{\rm\vskip2pt}
\newenvironment{remarknum}[1][]{\noindent\rm\textit{Remark\ #1}\ :\
}{\rm\vskip2pt}

\newenvironment{remarque}[1][]{\noindent\rm\textit{Remark}\ :\
}{\rm\vskip2pt}
\newenvironment{remarquenum}[1][]{\noindent\rm\textit{Remark\ #1}\ :\
}{\rm\vskip2pt}

\newenvironment{corollaire}[1][]{\noindent\sl\textbf{Corollary\ #1.}\
}{\rm\vskip3pt}

\begin{document}

\title[Harmonic functions on the real hyperbolic ball]{Harmonic functions 
on the real hyperbolic ball II\\ Hardy and Lipschitz spaces}
\author{Sandrine GRELLIER and Philippe JAMING}
\address{Universit\'e d'Orl\'eans\\ Facult\'e des Sciences\\ 
D\'epartement de Math\'ematiques\\BP 6759\\ F 45067 ORLEANS Cedex 2\\
FRANCE}
\email{grellier@labomath.univ-orleans.fr and
jaming@labomath.univ-orleans.fr}
\date{\today}
\thanks{The authors wish to thank A. bonami for valuable conversations
and advices.\\
Authors partially supported by the {\it European Commission}
(TMR 1998-2001 Network {\it Harmonic Analysis}).}
\begin{abstract} 
{\bf In this paper, we pursue the study of harmonic functions on the
real hyperbolic ball started in \cite{Ja3}. Our focus here is on the
theory of Hardy, Hardy-Sobolev and Lipschitz spaces of these functions.
We prove here that these spaces admit Fefferman-Stein like
characterizations in terms of maximal and square functionals. We further
prove that the hyperbolic harmonic extension of Lipschitz functions on
the boundary extend into Lipschitz functions on the whole ball.}

\end{abstract}

\maketitle

\section{Introduction}

In this article, the sequel of \cite{Ja3}, we study Hardy, Hardy-Sobolev
and Lipschitz spaces of harmonic functions on
the real hyperbolic ball. There are two main motivations for doing so\,:

While studying Hardy spaces of Euclidean harmonic functions on the unit
ball $\B_n$ of $\R^n$, one is often lead to consider estimates of these
functions on balls with radius smaller than the distance of the center
of that ball to the boundary $\S^{n-1}$ of $\B_n$. Thus hyperbolic
geometry is implicitly used in the study of Euclidean harmonic
functions.

The second motivation of this paper lies in the recent developments
of the theory of Hardy and Hardy-Sobolev
spaces of $\mm$-harmonic functions related to the complex hyperbolic
metric on the unit ball, as exposed in \cite{ABC} and \cite{BBG}. Our
aim here is to develop a similar theory in the case of the real
hyperbolic ball. In this paper, $n$ will be an integer, $n\geq3$ and
$p$ a real number, $0<p<+\infty$.

Our starting point is a result of \cite{Ja3} stating that the Hardy
spaces $\hh^p$ of hyperbolic harmonic functions ($\hh$-harmonic functions in the
terminology of \cite{Ja3}) admit an atomic decomposition similar to the
Euclidean harmonic functions. Then, for $0<p<+\infty$, define the space
$H^p(\S^{n-1})$ as $L^p(\S^{n-1})$ if $p>1$ and as the equivalent of
Garnett-Latter's atomic $H^p$-space if $0<p\leq1$ ({\it see} \cite{Ja3}
for the exact definition). This space has been characterized in terms of
square functionals of the Euclidean harmonic
extensions of its elements by Colzani \cite{Col}. We will here give 
these Fefferman-Stein characterizations directly in terms of their
$\hh$-harmonic extensions. More precisely, for an $\hh$-harmonic
function $u$, we prove the expected equivalence between $u\in\hh^p$
and its non-tangential maximal function, area integral or
Littlewood-Paley $g$-function belonging to $L^p(\S^{n-1})$.

In doing so, a choice of two methods is presented to us. We may either
use the link between $\hh$-harmonic functions and Euclidean harmonic
functions from \cite{Ja3} as for the atomic decomposition or else, adapt
the proofs in Fefferman-Stein \cite{FS1} to our context. In both cases
some difficulties appear. 

For the first method, the link we use only
allows to transfer results from the interior of the hyperbolic ball 
to the interior of the Euclidean ball, and from there to the boundary
$\S^{n-1}$ (by usual methods). Unfortunately a converse
link that would allow us to go back from the Euclidean ball to
the hyperbolic ball is only available in even dimension. Note also that
another link back from the Euclidean ball to the hyperbolic ball has been
exhibited in \cite{Sam} --see \cite{Ja3}, lemma 9-- but
this link implies loss of regularity and is thus not adapted to this
context.

\begin{figure}[ht]
\begin{center}
\includegraphics{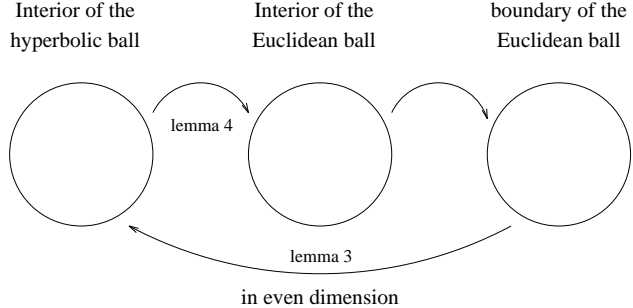}
\end{center}
\caption{Links between hyperbolic and Euclidean harmonic functions}
\end{figure}

In order to present unified proofs independent of the parity of the
dimension of the balls, we have thus inclined for Fefferman-Stein's
method. In doing so, the main difficulty is that the hyperbolic
Poisson kernels do not form a semi-group under convolution. In
particular, if the function $u$ is $\hh$-harmonic, the function
$u_\delta: x\mapsto u(\delta x)$ may not be $\hh$-harmonic anymore. 
This leads us to introduce the concept of $\hh_\delta$-harmonic functions and to
get estimates on these functions.

Our next interest is in developing a theory of Hardy-Sobolev spaces
of $\hh$-harmonic functions, similar to the one developed in
\cite{BBG}. The first step is to prove mean value inequalities for
$\hh$-harmonic functions and their derivatives. This is done by
adapting the proof in \cite{BBG} using the theory of hypo-elliptic
operators. We think that our mean value inequalities have an interest
in their own and that the proof should adapt to all rank one spaces of the
non-compact type. The remaining of the proofs are direct adaptations of
\cite{BBG}. However, as in \cite{Ja3} where it is proved that the 
boundary behavior of derivatives of $\hh$-harmonic functions dependeds 
on the parity of the dimension of $\B_n$, it is proved here that
the characterizations of Hardy-Sobolev spaces depend on the parity of
the order of derivation. Note that Graham \cite{Gra} has already notice
a dependance of the behavior of harmonic functions on the parity of 
the dimension of the balls.

Finally, we take advantage of the link between Euclidean and
hyperbolic harmonic functions to show how results on Lipschitz spaces
of Euclidean harmonic functions (see \cite{Gre}) can be transfered to the
hyperbolic harmonic context. In particular, we show that the
$\hh$-harmonic extension of a Lipschitz function on the boundary is
still a Lipschitz function of the same order on the whole ball. Further, 
we prove that in odd dimension, the limit-class preserved by $\hh$-harmonic
Poisson integrals is the Zygmund class of order $n$. This completes a result 
in \cite{Ja3} that states that this regularity is optimal in the sense that 
the $\hh$-harmonic extension of a function on $\S^{n-1}$ is at most in this
class.
\vskip6pt
This article is organized as follows. In the next section we present the
setting of our problem and state our main results. Section 3
is devoted to the proofs of the technical lemmas we will need, 
including the mean value inequalities. In section 4 we prove the 
Fefferman-Stein characterization of our $\hh^p$ spaces. The following
section is devoted to the proofs of similar characterizations for
Hardy-Sobolev spaces while in the last section we give
the results on Lipschitz spaces.
 
\section{Statement of the problem and results}

\subsection{$SO(n,1)$ and its action on $\B_n$}
Let $G=SO(n,1)\subset GL_{n+1}(\R)$, ($n\geq3$) be the identity
component of the group of matrices
$g=(g_{ij})_{0\leq i,j\leq n}$ such that $g_{00}\geq1$, $\det g=1$ and
that leave invariant the quadratic form $-x_0^2+x_1^2+\ldots+x_n^2$.
$G$ admits a Cartan decomposition $G=K\overline{A_+}K$ where 

$$K=\left\{k=\left(\begin{matrix} 1&0\\0&\hat k\\ \end{matrix}\right)\ :\ \hat
k\in SO(n)\right\}$$
and

$$A_+=\left\{a_t=\left(\begin{matrix}\cosh t&\sinh t&0\\ \sinh t&\cosh
t&0\\ 0&0&Id_{n-1}\\ \end{matrix}\right)\ :\ t\in\R_+\right\}.$$
In this decomposition, every $g\in G$ can be written
$g=k_ga_{t(g)}k^{'}_g$.

Let $\abs{.}$ be the Euclidean norm on $\R^n$ and $\scal{.,.}$ the
associated scalar product. Let $\B_n=\{x\in\R^n\
:\ \abs{x}<1\}$ and $\S^{n-1}=\partial\B_n=\{x\in\R^n\ :\ \abs{x}=1\}$.
The homogeneous space $G/K$ can be identified with $\B_n$, and it is well known 
({\it see} \cite{Sam}) that $SO(n,1)$ acts conformaly on $\B_n$
by $y=g.x$ with

$$y_p=\frac{\frac{1+\abs{x}^2}{2}g_{p0}+\sum_{l=1}^ng_{pl}x_l}{\frac{1-\abs{x}^2}{2}+
\frac{1+\abs{x}^2}{2}g_{00}+\sum_{l=1}^ng_{0l}x_l}
\qquad\mathrm{for\ }p=1,\ldots,n.$$
The invariant measure on $\B_n$ is given by

$$d\mu=\frac{dx}{(1-\abs{x}^2)^{n-1}}=
\frac{r^{n-1}drd\sigma}{(1-r^2)^{n-1}}$$
where $dx$ is the Lebesgue measure on $\B_n$ and $d\sigma$ is the
surface measure on $\S^{n-1}$.

We will need the following elementary facts about this action (see \cite{Jathese}):

{\bf Fact 1.} {\sl Let $g\in SO(n,1)$ and let $x_0=g.0$. If
$0<\eps<\frac{1}{6}$, then}

$$B\bigl(x_0,\frac{\sqrt{2}}{8}(1-\abs{x_0}^2)\eps\bigr)\subset
g.B(0,\eps)\subset B\bigl(x_0,6(1-\abs{x_0}^2)\eps\bigr).$$

{\bf Fact 2.} {\sl Let $g\in SO(n,1)$ and let $x_0=g.0$. Let $v$ be a
smooth function on $\B_n$ and  define $f$ on $\B_n$ by $f(x)=v(g.x)$.
Then, for every $k$,

$$(1-\abs{x_0}^2)^k\abs{\nabla^kv(x_0)}\leq C\abs{\nabla^kf(0)},$$
where $\abs{\nabla^k}$ means 
$\sup\{\abs{\frac{\partial^{\abs{\alpha}}}{\partial x^\alpha}}\ :\
\abs{\alpha}\leq k\}$.}
\subsection{The invariant laplacian on $\B_n$ and the associated Poisson
kernel}

From \cite{Sam} (see also \cite{Ey1},\cite{Ey2}), we know that the invariant laplacian on $\B_n$
for the considered action can be written as

$$D=(1-r^2)^2\Delta+2(n-2)(1-r^2)\sum_{i=1}^nx_i\frac{\partial}{\partial
x_i}$$
where $r=\abs{x}=(x_1^2+\ldots+x_n^2)^{1/2}$ and $\Delta$
is the Euclidean laplacian
$\displaystyle\Delta=\sum_{i=1}^n\frac{\partial^2}{\partial x_i^2}$.

Note that $D$ is given in radial-tangential coordinates by

$$D=\frac{1-r^2}{r^2}\ent{(1-r^2)N^2+(n-2)(1+r^2)N+(1-r^2)\Delta_\sigma}$$
with $\displaystyle N=r\frac{d}{dr}=\sum_{i=1}^nx_i\frac{\partial}{\partial x_i}$ and
$\Delta_\sigma$ the tangential part of the
Euclidean laplacian.

\begin{definition} A function $u$ on $\B_n$ is
$\hh$-harmonic if $Du=0$ on $\B_n$.
\end{definition}

\begin{notation} Let
$L=\frac{1}{r^2}\ent{(1-r^2)N^2+(n-2)(1+r^2)N+(1-r^2)\Delta_\sigma}$. 
Thus $Du=0$ if and only if $Lu=0$.
\end{notation}

Green's formula for $D$ is given by the following theorem :

\begin{theoremnom}[1]{Green's formula} Let $\Omega$ be an open
subset of
$\B_n$ with $\cc^1$ smooth boundary and let $\vec{n}$ be the exterior
normal
to $\partial\Omega$. Then for every functions $u,v\in \cc^2(\Omega)$,

\begin{align}
\int_\Omega(1-\abs{x}^2)^{-n}(uDv-vDu)dx=&
\int_{\partial \Omega}\ent{u\frac{\partial v}{\partial
\vec{n}}-v\frac{\partial u}{\partial\vec{n}}}(1-r^2)^{-n+2}
d\sigma.\notag
\end{align}
\end{theoremnom}

The Poisson kernel that solves the Dirichlet problem  associated to $D$
is given by

$$\P_h(r\eta,\xi)=\left(\frac{1-r^2}{1+r^2-2r\scal{\eta,\xi}}
\right)^{n-1}$$
for $0\leq r<1$, $\eta,\xi\in\S^{n-1}$ {\it i.e.} for
$r\eta\in\B_n$ and $\xi\in\S^{n-1}$.

Recall that the Euclidean Poisson kernel on the ball is given by

$$\P_e(r\eta,\xi)=\frac{1-r^2}{(1+r^2-2r\scal{\eta,\xi}
)^{\frac{n}{2}}}$$

\begin{notation} For a distribution $\ffi$ on $\S^{n-1}$, we define
$\P_e\ent{\ffi}:\B_n\mapsto\R$ and $\P_h\ent{\ffi}:\B_n\mapsto\R$ by

\begin{align}
\P_e\ent{\ffi}(r\eta)=&\scal{\ffi,\P_e(r\eta,.)}\notag\\
\P_h\ent{\ffi}(r\eta)=&\scal{\ffi,\P_h(r\eta,.)}\notag
\end{align}
$\P_e\ent{\ffi}$ is the {\it Poisson integral} of $\ffi$, and
$\P_h\ent{\ffi}$ will be called the {\it $\hh$-Poisson integral} of
$\ffi$.
\end{notation}

\subsection{Expansion of $\hh$-harmonic functions in spherical
harmonics}

Let ${}_2F_1$ denote Gauss' {\it
hyper-geometric function} and
let $F_l(x)={}_2F_1(l,1-\frac{n}{2},l+\frac{n}{2};x)$ and
$f_l(x)=\frac{F_l(x)}{F_l(1)}$. (See \cite{Erd} for properties of
${}_2F_1$ used here).

In \cite{Mi1}, \cite{Mi2} and \cite{Sam}, the spherical harmonic
expansion of $\hh$-harmonic functions has been obtained. Another 
proof, based on the method developped in \cite{ABC} for $\mm$-harmonic
functions, can be found in \cite{Jathese}. We have the 
following :

\begin{theorem}[2] Let $u$ be an $\hh$-harmonic function of class
$\cc^2$ on $\B_n$. Then the spherical harmonic expansion of $u$ is given
by

$$u(r\zeta)=\sum_lf_l(r^2)u_l(r\zeta),$$
where this series is absolutely convergent and uniformaly convergent
on every compact subset of $\B_n$.
\end{theorem}

Moreover, if we denote by $\Z_l^\zeta$ the zonal function of
order $l$ with pole $\zeta$, then the hyperbolic Poisson kernel is given
by

$$\P_h(r\zeta,\xi)=\sum_{l\geq0}\frac{F_l(r^2)}{F_l(1)}
r^l\Z_l^\zeta(\xi).$$
Recall also that the Euclidean Poisson kernel is given by

$$\P_e(r\zeta,\xi)=\sum_{l\geq0}r^l\Z_l^\zeta(\xi).$$
In case the dimension $n$ is {\it even}, this two kernels are linked by
the following.

\begin{lemmanum}[3] Assume $n$ is even, and write $n=2p$. There 
exists $p$ polynomials $P_0,P_1,\ldots,P_{p-1}$ such that, for
every $r\zeta\in\B_n,\xi\in\S^{n-1}$,

$$\P_h(r\zeta,\xi)=\sum_{k=0}^{p-1}P_k(r)(1-r^2)^k
\frac{\partial^k}{\partial r^k}\P_e(r\zeta,\xi).$$
\end{lemmanum}

\begin{proof}[Proof] For $a\in\R$, write
$(a)_k=\frac{\Gamma(a+k)}{\Gamma(a)}$. From \cite{Erd} we get

$$F_l(x)={}_2F_1(l,1-p,l+p,x)=\frac{1}{(l+p)_{p-1}}
\frac{(1-x)^{2p-1}}{x^{l+p-1}}
\frac{d^{p-1}}{dx^{p-1}}\bigl(x^{l+2(p-1)}(1-x)^{-p}\bigr).$$

Let $\alpha_{l,j}$ be defined by $\alpha_l,0=1$ and
$\alpha_{l,j+1}=\bigl(l+2(p-1)-j\bigr)\alpha_{l,j}$, then by Leibniz'
formula

$${}_2F_1(l,1-p,l+p,x)=\frac{1}{(l+p)_{p-1}}\sum_{j=0}^{p-1}
\begin{pmatrix} p-1\\ j\\ \end{pmatrix}
(p)_j\alpha_{l,j}x^j(1-x)^j.$$
In particular ${}_2F_1(l,1-p,l+p,1)=\frac{1}{(l+p)_{p-1}}$ thus

$$\frac{F_l(x)}{F_l(1)}=\sum_{j=0}^{p-1}
\begin{pmatrix}p-1\\ j\\ \end{pmatrix}(p)_j\alpha_{l,j}x^j(1-x)^j.$$
Furthermore, it is easy to see that one can write

$$\alpha_{l,j}=\sum_{k=0}^ja_{k,j}l(l-1)\ldots(l-k+1)$$
where the coefficients $a_{k,j}$ are independent from $l$. It results
from this and the spherical harmonics expansions of $\P_h$ and $\P_e$
that there exist polynomials $P_0,P_1,\ldots,P_{p-1}$ such that, for
every $r\zeta\in\B_n,\xi\in\S^{n-1}$,

$$\P_h(r\zeta,\xi)=\sum_{k=0}^{p-1}P_k(r)(1-r^2)^k
\frac{\partial^k}{\partial r^k}\P_e(r\zeta,\xi).$$
which completes the proof.\hfill$\Box$
\end{proof}

In \cite{Ja3}, the following link between euclidean harmonic functions
and $\hh$-harmonic functions has been exhibited :

\begin{lemmanum}[4] There exists a function
$\eta:\ent{0,1}\times\ent{0,1}\mapsto\R^+$ such that
\begin{description}
\item[i] $\P_e(r\zeta,\xi)=\int_0^1\eta(r,\rho)\P_h(\rho
r\zeta,\xi)d\rho$,

\item[ii] for every $k$, there exists a constant $C_k$ such that for every
$r\in\ent{0,1}$, 

$$\int_0^1\abs{\left(r\frac{\partial}{\partial 
r}\right)^k\eta(r,\rho)}d\rho\leq \frac{C}{(1-r)^k}.$$
\end{description}
\end{lemmanum}

\begin{proof}[Proof] According to \cite{Ja3}, the function $\eta$ is 
given by

$$\eta(r,s)=c(1-r^2)(1-r^2s^2)^{2-n}\ent{(1-s)(1-sr^2)}^{\frac{n}{2}-2}
s^{\frac{n}{2}-1}.$$
The etimate {\it ii/} is easily obtained by differentiation.\hfill$\Box$
\end{proof}

\subsection{Hardy and Hardy-Sobolev spaces}

The aim of this article is to extend Fefferman-Stein \cite{FS1} theory
to Hardy and Hardy-Sobolev spaces of $\hh$-harmonic functions. We will
therefore need to define analogs of non-tangential maximal
functions, area integrals and Littlewood-Paley $g$ functions.

\begin{definitionfr} For $0<\alpha<1$ and $\zeta\in\S^{n-1}$, let 
$\aa_\alpha(\zeta)$ be the interior of the convex hull of $B(0,\alpha)$
and
$\zeta$ ; $\aa_\alpha(\zeta)$ will be called {\it non-tangential
approach region}.
\begin{figure}[ht]
\begin{center}
\includegraphics{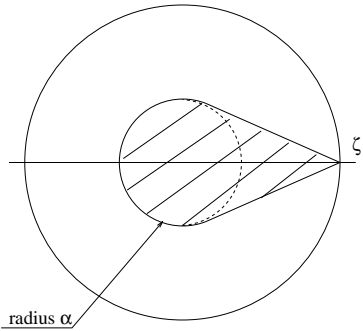}
\end{center}
\caption{non-tangential approach region $\aa_\alpha(\zeta)$}
\end{figure}

For a function $u$ defined on $\B_n$, define the following functions on
$\S^{n-1}$ :

\begin{enumerate}
\item $\mm\ent{u}(\xi)=\sup_{0<r<1}\abs{u(r\xi)}$,

\item $\mm_\alpha\ent{u}(\xi)=\sup_{x\in\aa_\alpha(\xi)}\abs{u(x)}$.

\item $\displaystyle S_\alpha\ent{u}(\xi)=\ent{\int_{\aa_\alpha(\xi)}
\abs{\nabla u(x)}^2(1-\abs{x}^2)^{-n+2}dx}^{\frac{1}{2}}$.

\item $\displaystyle S_\alpha^N\ent{u}(\xi)=\ent{\int_{\aa_\alpha(\xi)}
\abs{N u(x)}^2(1-\abs{x}^2)^{-n+2}dx}^{\frac{1}{2}}$.

\item $g\ent{u}(\xi)=\ent{\int_0^1\abs{\nabla
u(t\xi)}(1-t^2)dt}^{\frac{1}{2}}$.

\item $g^N\ent{u}(\xi)=\ent{\int_0^1\abs{N u(t\xi)}(1-t^2)
dt}^{\frac{1}{2}}$.
\end{enumerate}
\end{definitionfr}
We can then define the Hardy spaces for $0<p<+\infty$ as

$$\hh^p=\{u\ \mathrm{\hh-harmonic}\ :\ \mm\ent{u}\in L^p(\S^{n-1})\}.$$

We will prove the following result :

\begin{theorem}[A] For $0<p<2$ and $u$ $\hh$-harmonic, the
following are equivalent :

\begin{enumerate}
\item $u\in\hh^p$.

\item $u$ has a boundary distribution in $H^p(\S^{n-1})$.

\item $\mm_\alpha\ent{u}\in L^p(\S^{n-1})$ for some $0<\alpha<1$.

\item $S_\alpha\ent{u}\in L^p(\S^{n-1})$ for some $0<\alpha<1$.

\item $S_\alpha^N\ent{u}\in L^p(\S^{n-1})$ for some $0<\alpha<1$.

\item $g\ent{u}\in L^p(\S^{n-1})$,

\item $g^N\ent{u}\in L^p(\S^{n-1})$.
\end{enumerate}
Moreover, the equivalence of $1$, $2$ and $3$ is valid for $0<p<+\infty$.
\end{theorem}

\begin{remark} This theorem implies, in particular, that if assertions
$3$, $4$ and $5$ are satisfied for some $\alpha$, they are satified for
every $\alpha$.
\end{remark}

Note that with lemma 3, part of this theorem is obvious 
in case the dimension $n$ is even. However, we prefer giving 
here unified proofs independent of the parity of the dimension.

Define now the Hardy-Sobolev spaces for $0<p<+\infty$ and $k\in\N$ as

$$\hh^p_k=\{u\ \hh-\mathrm{harmonic}\ :\ \mathrm{for\ all}
\ j\leq k,\ \mm\ent{\nabla^ju}\in L^p(\S^{n-1})\}.$$
and

$$H^p_k(\S^{n-1})=\{f\in H^p(\S^{n-1})\ ;\ \nabla^jf\in
H^p(\S^{n-1}),\ 0\leq j\leq k\}.$$

We prove the following theorem :

\begin{theorem}[B] For $0<p<2$, for every integer 
$0\leq k\leq n-2$ and for every $\hh$-harmonic function $u$, 
the following are equivalent :

\begin{enumerate}
\item $u\in\hh^p_k$.

\item $u$ has a boundary distribution in $H^p_k(\S^{n-1})$.

\item $u$ has a boundary distribution $f$ satisfying
$(-\Delta_\sigma)^{\frac{l}{2}}f\in H^p(\S^{n-1})$ for $0\leq l\leq k$.

\item $u\in\hh^p_{k-1}$ and for some $\alpha$ such that $0<\alpha<1$,
$\mm_\alpha\ent{(-\Delta_\sigma)^{k/2}u}\in L^p(\S^{n-1})$.

\item $u\in\hh^p_{k-1}$ and for some $\alpha$ such that $0<\alpha<1$,
$S_\alpha\ent{(-\Delta_\sigma)^{k/2}u}\in L^p(\S^{n-1})$.

\item $u\in\hh^p_{k-1}$ and for some $\alpha$ such that $0<\alpha<1$,
$S_\alpha^N\ent{(-\Delta_\sigma)^{k/2}u}\in L^p(\S^{n-1})$.

\item $u\in\hh^p_{k-1}$ and for some $\alpha$ such that $0<\alpha<1$,
$S_\alpha\ent{N^ku}\in L^p(\S^{n-1})$.

\item $u\in\hh^p_{k-1}$ and for some $\alpha$ such that $0<\alpha<1$,
$S_\alpha^N\ent{N^ku}\in L^p(\S^{n-1})$.

\item $u\in\hh^p_{k-1}$ and for some $\alpha$ such that $0<\alpha<1$,
$S_\alpha\ent{\nabla^ku}\in L^p(\S^{n-1})$.
\end{enumerate}

Moreover, the equivalence of $1$, $2$, $3$ and $4$ is valid for $0<p<+\infty$.
\end{theorem}

\begin{remarknum}[1] Again, this theorem implies that if assertions 
$4$ to $9$ are satisfied for some $\alpha$, they are satified for
every $\alpha$.
\end{remarknum}

\begin{remarknum}[2] As $(-\Delta_\sigma)^{1/2}$ preserves
$\hh$-harmonicity, the equivalence of $2$, $4$, $5$ and $6$ means that 
$(-\Delta_\sigma)^{1/2}u\in\hh^p$. The equivalence between $2$ and $3$
then follows from the atomic decomposition of $H^p(\S^{n-1})$ and 
standard singular integral arguments.
\end{remarknum}

\begin{remarknum}[3] Let $\lll_{i,j}=x_i\frac{\partial}{\partial x_j}
-x_j\frac{\partial}{\partial x_i}$. They satisfy
$\Delta_\sigma=\sum_{i<j}\lll_{i,j}^2$. It is obvious that the
$\lll_{i,j}$'s commute with the invariant laplacian $D$ so that they
preserve $\hh$-harmonicity. Further, if $l$ is an odd integer with 
$l=2l_0+1$, we can replace $(-\Delta_\sigma)^{l/2}$ in $4$, $5$ and $6$ by the 
set of operators 
$\{\Delta_\sigma^{l_0}\lll_{i,j}u\ :\ 0\leq i,j\leq n\}$.
\end{remarknum}

Define now

$$\hh^p_{k,N}=\{u\in\hh^p;\ \mm\ent{N^lu}\in L^p(\S^{n-1}),
\ 0\leq l\leq k\}.$$
We study the relationship between $\hh^p_k$ and $\hh^p_{k,N}$. The
situation is slightly different as parity of the order of derivation
is involved.

\begin{theorem}[C] For $0<\alpha<1$, $0<p<+\infty$, and $k$ an integer,
$0\leq k\leq n-2$. Then
\begin{enumerate}
\item If $k$ is {\rm even}, the following are equivalent :

\begin{enumerate}
\item $u\in\hh^p_{k,N}$.

\item For some  $\alpha$ such that $0<\alpha<1$, for every $0\leq l\leq k$,
$\mm_\alpha\ent{N^l u}\in L^p(\S^{n-1})$.

\item $u\in\hh^p_k$ and hence all the equivalent properties stated in
theorem B are valid.
\end{enumerate}

\item If $k$ is {\rm odd}, the following are equivalent :

\begin{enumerate}
\item $u\in\hh^p_{k,N}$

\item For some  $\alpha$ such that $0<\alpha<1$, for every $0\leq l\leq k$,
$\mm_\alpha\ent{N^l u}\in L^p(\S^{n-1})$.

\item $u\in\hh^p_{k-1}$ and 
$\mm\ent{(1-r^2)\Delta_\sigma^{\frac{k+1}{2}}u}\in L^p(\S^{n-1})$.

\end{enumerate}
\end{enumerate}
\end{theorem}

\begin{remark} The third assertion in part 2, {\it i.e.} when $k$ is
odd,  is in particular satisfied
when $u\in\hh^p_{k-1}$ and
$\mm_\alpha\ent{(-\Delta_\sigma)^{\frac{k}{2}}u}\in 
L^p(\S^{n-1})$, that is when $u$ is in $\hh^p_k$.
This gives the inclusion 
$\hh^p_k\subset\hh^p_{k,N}$. However, the space $\hh_{k,N}^p$
is strictly bigger. 

This result looks, at first sight, quite surprising
since it is usually expected that the radial derivative dominates the
gradient. In fact, this is naturally true in the interior of the domains
and for instance, $S_\alpha(N^lu)\in L^p(\S^{n-1})$, $0\leq l\leq k$
implies (and in fact is equivalent to) $u\in\hh^p_k$ as stated in
theorem B.

It is no longer true for conditions involving the behavior of the radial
derivatives near the boundary. For instance, when $k=1$, we know
from \cite{Ja3} that $Nu$ has a boundary distribution that is 
identically zero. So, for $u$ $\hh$-harmonic, to be in $\hh^p_{1,N}$
can not be translated as a constraint on the boundary behaviour of $u$.
\end{remark}

\section{Preliminary lemmas}

\subsection{Mean value inequalities}
\label{inmoy}

Recall that $\hh$-harmonic functions satisfy the following mean value
equalities :

{\sl Let $a\in\B_n$ and $g\in SO(n,1)$ such that $g.0=a$. Then, for 
every $\hh$-harmonic function $u$,}

$$u(a)=\frac{1}{\mu\bigl(B(0,r)\bigr)}
\int_{g.B(0,r)}u(x)d\mu(x).$$
Thus, with fact 1 and $d\mu=\frac{dx}{1-\abs{x}^2}$, we get

$$\abs{u(a)}\leq\frac{C}{(1-\abs{a}^2)^n}
\int_{B\bigl(a,6(1-\abs{a}^2)\eps\bigr)}\abs{u(x)}dx
\leqno(3.1)$$

We will also need mean
value inequalities for normal derivatives of $\hh$-harmonic functions, 
in particular when we study Hardy-Sobolev spaces. But, normal
derivatives
of $\hh$-harmonic functions are no longer $\hh$-harmonic, so that
inequality (3.1) does not apply to them. 

To obtain this inequalities, we will 
follow the main lines of the proof in \cite{BBG} for $\mm$-harmonic
functions.

Therefore, we will first study the commutator between $N^k$ and $L$
(which is easier to compute than the commutator between $N^k$ and $D$).
This leads us to the existence of an elliptic operator $\N_q$ such
that
for every $\hh$-harmonic function $u$, $N^ku$ is annihilated by $\N_q$.
We can then apply $L^2$ theory of elliptic operators and get
estimates 
for $N^ku$ in $0$. To obtain the estimates in an arbitrary point $a$ of
$\B_n$, we transport the result from $0$ to $a$ with help of the action
of 
$SO(n,1)$ on $\B_n$ by computing the action of $g\in SO(n,1)$ on $\N_q$.

Note that

$$LN-NL=2L+2(N^2+\Delta_\sigma)-2(n-2)N\leqno(3.2)$$
Moreover an easy induction argument shows that there exist
two sequences of polynomials
$(P_k)_{k\geq1}$ and $(Q_k)_{k\geq1}$ of degree $k-1$ such that for
$k\geq1$,

$$LN^k=(N+2I)^kL+P_k(N)N^2+Q_k(N)\Delta_\sigma-2(n-2)(N+2I)^{k-1}N.$$

From this, using the same induction as in \cite{BBG}, we get

\begin{proposition}[5] For every $k$, there exist polynomials
$S_k(x,y)$
of degree at most $q-1$ (with $q=2^{k-1}$) and $R_k(x,y)=x^q+\ldots$
such that,
if $u$ is $\hh$-harmonic, then

$$L\bigl(R_k(L,\Delta_\sigma)-S_k(L,\Delta_\sigma)N\bigr)N^k u=0.$$
\end{proposition}

We thus conclude that if $u$ is $\hh$-harmonic, then $v=N^ku$
is a solution of an equation $\N_qv=0$ with $\N_q=L^{q+1}+\ldots$ and
$q=2^{k-1}$.

We will use the following formalism :
if $\M$ is a differential operator and $\Phi$ a diffeomorphism of $\B_n$
and 
if $f=v\circ\Phi$, then define $\Phi^*\M$ by

$$\Phi^*\M(f)=(\M v)\circ\Phi.$$
It is then obvious that

\begin{align}
\Phi^*(\M_1\circ M_2)=&(\Phi^*\M_1)\circ(\Phi^*\M_2)
\tag{3.3}\\
\Phi^*(h\M)=&h(\Phi).\Phi^*\M\tag{3.4}
\end{align}

We will consider $v=N^ku$ with $u$ $\hh$-harmonic, so that $\N_qv=0$
with
$\N_q=L(R_k(L,\Delta_\sigma)-P_k(L,\Delta_\sigma)N)=L^{q+1}+\ldots$

Let $g\in SO(n,1)$ be such that $g.0=\rho\zeta=a\in\B_n$ and let
$\Phi_a:\matrix \B_n&\mapsto&\B_n\cr x&\mapsto&g.x\cr\endmatrix$\ .

But, by definition, $D$ is invariant by the action of $SO(n,1)$ on 
$\B_n$, that is $\Phi_a\!{}^*D=D$. On the other hand, $D=(1-\abs{x}^2)L$
thus (3.4) tells us that 
$\Phi_a\!{}^*D=(1-\abs{\Phi_a(x)}^2)\Phi_a\!{}^*L$,
which implies that $\Phi_a\!{}^*L=\frac{1-\abs{x}^2}{1-\abs{g.x}^2}L$,
and the formula of \cite{Sam} page 39 gives

$$\Phi_a\!{}^*L=\frac{(1+\rho^2\abs{x}^2-2\rho\scal{x,\zeta}
)^2}{1-\rho^2}L.$$

Further $\Phi_a\!{}^*N$ is a differential operator of order $1$ with
$\cc^\infty$ coefficients defined by

\begin{align}
\Phi_a\!{}^*Nf(x)=&\scal{\Phi_a(x),dv_{\Phi_a(x)}}=
\scal{\Phi_a(x),d(f\circ\Phi_a^{-1})_{\Phi_a(x)}}\notag\\
=&\scal{\Phi_a(x),df_z.d(\Phi_a^{-1})_{\Phi_a(x)}}\notag
\end{align}
thus, if $x\in B(0,\eps)$ then, with fact 1,
$\Phi_a(x)\in B\bigl(a,6(1-a^2)\eps\bigr)$, and with fact 2 (applied to
$v(x)=x$),
the coefficients of $(1-\abs{a}^2)\Phi_a\!{}^*N$ as well as 
their derivatives are  $\cc^\infty$ and bounded independently of $a$.

As $\Phi_a\!{}^*N^2=(\Phi_a\!{}^*N)\circ(\Phi_a\!{}^*N)$,
$(1-\abs{a}^2)^2\Phi_a\!{}^*N^2$ is a differential operator of order $2$
with $\cc^\infty$ coefficients bounded (as well as their derivatives)
independently of $a$.

At last,
$\Delta_\sigma=\frac{r^2}{1-r^2}L-N^2-(n-2)\frac{1+r^2}{1-r^2}N$
thus $(1-\abs{a}^2)\Phi_a\!{}^*\Delta_\sigma$ is also a differential
operator
of order $2$ with $\cc^\infty$ coefficients bounded (as well as their 
derivatives) independently of $a$.

Finally, $\N_q=L^{q+1}+$ terms of order $\leq 2q$ in $L,\Delta_\sigma$
and $N$ with  $\cc^\infty$ coefficients, thus

\begin{align}
\Phi_a\!{}^*\N_q=&\Phi_a\!{}^*L^{q+1}+\text{terms of order }\leq2q\text{
with }
\cc^\infty\text{ coefficients}\notag\\
=&\frac{(1+\rho^2\abs{x}^2-2\rho\scal{x,\zeta})^{2(q+1)}}{(1-\rho^2)^{q+1}}L^{q+1}
\notag\\
&+\text{terms of order }\leq2q\text{ with }
\cc^\infty\text{ coefficients}\notag\\
=&\frac{\abs{\rho x-\zeta}^{2(q+1)}}{(1-\rho^2)^{q+1}}L^{q+1}+\R_{q,a}
\notag\end{align}
where $\R_{q,a}$ is a differential operator of order $\leq2q$ with
$\cc^\infty$
coefficients.

Let $u$ be an $\hh$-harmonic function, $v=N^ku$ and $f=v\circ\Phi_a$. As
$v$ satisfies $\N_qv=0$, $f$ satisfies
$(1-\abs{a}^2)^{q+1}\Phi_a\!{}^*\N_qf=0$ on $B(0,\eps)$ (with 
{\it e.g.} $\eps<\frac{1}{6}$). We have thus shown that
$(1-\rho^2)^{q+1}\Phi_a\!{}^*\N_q$ satisfies on $B(0,\eps)$,
$\eps<1/6$, all the hypotheses (with constants independent on $a$)
of the following theorem (see \cite{BBG} page 678)\,:

\begin{theorem}[6] Suppose $P(d)$ is a differential operator in $\R^N$,

$$P(D)=\sum_{\abs{\alpha}\leq2q}h_\alpha(x)D^\alpha\quad\mathrm{where}\quad
D^\alpha=\frac{\partial^{\alpha_1}}{\partial x_1^{\alpha_1}}\ldots
\frac{\partial^{\alpha_N}}{\partial x_N^{\alpha_N}},$$
which is elliptic with constant $c_0$ in $B(0,\eps)$, that is,

$$\sum_{\abs{\alpha}=2q}h_\alpha(x)\xi^\alpha\geq c_0\abs{\xi}^{2q}\quad\mathrm{for}\quad
\xi\in\R^N,$$
and with $h_\alpha\in\cc^{\infty}\bigl(\overline{B(0,\eps)}\bigr)$. Assume that $P(D)f=0$ in $B(0,\eps)$.
Then, for all non-negative integers $m$ and all $p$ such that $0<p<\infty$,

$$\abs{\nabla^mf(0)}\leq C\left(\int_{\abs{x}\leq\eps}\abs{f(x)}^pdx\right)^{1/p},$$
where $C$ depends only on $c_0,\eps,m,p$ and a bound of the norms of the functions $h_\alpha$ in some
$\cc^l\bigl(\overline{B(0,\eps)}\bigr)$-space with $l=l(m)$
\end{theorem}

From this, we get

$$\abs{\nabla^df(0)}\leq
c\left(\int_{\abs{x}\leq\eps}\abs{f(x)}^pdx\right)^{\frac{1}{p}}\leq
c\left(\int_{\abs{x}\leq\eps}\abs{f(x)}^p\frac{dx}{(1-\abs{x}^2)^n}\right)
^{\frac{1}{p}}$$
or, with {\sl Fact 2},

$$\abs{\nabla^dv(a)}\leq
c\left(\int_{B(0,\eps)}\abs{v\circ\Phi_a(x)}^pd\mu(x)\right)^{\frac{1}{p}}
\times(1-\abs{a}^2)^{-d}$$
where $\mu$ is the $G$-invariant measure on $\B_n$. Thus

$$\abs{\nabla^dv(a)}\leq
c\left(\int_{g.B(0,\eps)}\abs{v(x)}^pd\mu(x)\right)^{\frac{1}{p}}
\times(1-\abs{a}^2)^{-d}$$
and, with {\sl Fact 1},

\begin{align}
\abs{\nabla^dv(a)}\leq&
c\left(\int_{B\bigl(a,6(1-\abs{a}^2)\eps\bigr)}\abs{v(x)}^p\frac{dx}{(1-\abs{x}^2)^n}
\right)^{\frac{1}{p}}\times(1-\abs{a}^2)^{-d}\notag\\
\leq&c(1-\abs{a}^2)^{-d-\frac{n}{p}}
\left(\int_{B\bigl(a,6(1-\abs{a}^2)\eps\bigr)}\abs{v(x)}^pdx\right)^{\frac{1}{p}}\notag
\end{align}
In conclusion, we have just proved the following lemma

\begin{lemmenom}[7]{Mean Value Inequality} For every
$0<\eps<\frac{1}{6}$, $k,d\in\N$, $0<p<+\infty$,
there exists a constant $c$ such that, for every $\hh$-harmonic function
$u$,
and every $a\in\B_n$,

$$\abs{\nabla^dN^ku(a)}\leq c(1-\abs{a})^{-d-\frac{n}{p}}\left(
\int_{B\bigl(a,6(1-\abs{a}^2)\eps\bigr)}\abs{N^ku(x)}^pdx\right)^{\frac{1}{p}}.$$
\end{lemmenom}

\begin{remarque} In case $d=0$ ($\nabla^0=I$), $k=0$ and $p=1$, 
we again obtain inequality (3.1).
\end{remarque}

Let $\lll_{i,j}=x_i\frac{\partial}{\partial
x_j}-x_j\frac{\partial}{\partial x_i}$ ($1\leq i<j\leq n$). One easily
sees 
that $N\lll_{i,j}=\lll_{i,j}N$ and that
$\Delta\lll_{i,j}=\lll_{i,j}\Delta$, 
thus  $L\lll_{i,j}=\lll_{i,j}L$ and $D\lll_{i,j}=\lll_{i,j}D$. In
particular,
if $u$ is $\hh$-harmonic, so is $\lll_{i,j}^ku$ for every $k\in\N$. 
Applying lemma 7 to $\lll_{i,j}u$, for every
$0<\eps<\frac{1}{6}$ 
and every $0<p<+\infty$, there exists a constant $C$ such that for every
$\hh$-harmonic function $u$, for every $1\leq i<j\leq n$ and every
$k\in\N$,
for every $d$, and every $a\in\B_n$,

$$\abs{\nabla^d\lll_{i,j}^ku(a)}\leq C(1-\abs{a})^{-d-\frac{n}{p}}\left(
\int_{B\bigl(a,6(1-\abs{a}^2)\eps\bigr)}\abs{\lll_{i,j}^ku(x)}^p
dx\right)^{\frac{1}{p}}.\leqno(3.5)$$

\begin{remarque} Let $\tilde\nabla^ku$ be defined by

$$\{\X N^qu\tq\X=\prod_{l=1}^p\lll_{i_l,j_l},\ p+q\leq k\},$$
then (3.5) implies that lemma 7 stays true if we
replace $N^k$ by $\tilde\nabla^ku$. But, outside a fixed neighborhood 
$V$ of $0$, $\abs{\tilde\nabla^ku}\simeq\abs{\nabla^ku}$, thus for every 
$a\in\B_n\setminus V$

$$\abs{\nabla^d\nabla^ku(a)}\leq C(1-\abs{a})^{-d-\frac{n}{p}}
\left(\int_{B\bigl(a,2(1-\abs{a}^2)\eps\bigr)}\abs{\nabla^ku(x)}^p
dx\right)^{\frac{1}{p}}.$$
As for $a\in V$ one can apply theorem 6 on
$B\bigl(a,2(1-\abs{a}^2)\eps\bigr)$ with constants independent of $a$,
we get the previous inequality on $V$ (recall that $\nabla^k$ means the
set of all derivatves of order less than $k$). We thus get the 
following proposition :
\end{remarque}

\begin{proposition}[8] For every $0<\eps<\frac{1}{6}$ and every
$0<p<+\infty$, there exists a constant $C$ such that for every
$\hh$-harmonic function $u$, every  $k\in\N$,  $d\geq0$, and for every
$a\in\B_n$,

$$\abs{\nabla^{k+d}u(a)}\leq C(1-\abs{a})^{-d-\frac{n}{p}}\left(
\int_{B\bigl(a,6(1-\abs{a}^2)\eps\bigr)}\abs{\nabla^ku(x)}^pdx\right)^{\frac{1}{p}}.$$
\end{proposition}

\begin{remarknum}[1] In the sequel, we will not distinguish anymore between
$\nabla^k$ and $\tilde\nabla^k$.
\end{remarknum} 

\begin{remarknum}[2] The previous inequality can be restated in an invariant
form under the action of the group, using invariant gradient and, more
generally covariant derivation.
%
if
%
\end{remarknum}

\begin{corollaire}[9] Let $0<\alpha<\beta<1$, $k,d\in\N$. Then
there
exists a constant $c$ such that for every $\hh$-harmonic function $u$,

$$\mm_\alpha\bigl((1-\abs{z})^d\nabla^dN^ku\bigr)\leq
c\mm_\beta(N^ku).$$
\end{corollaire}

\begin{proof}[Proof] This is an immediate consequence of lemma 7
and the fact that if $\alpha<\beta$ and if $\eps$ is small enough then,
for every $\xi\in\S^{n-1}$ and every $a\in\aa_\alpha(\xi)$,
$B(a,6(1-\abs{a}^2)\eps)\subset\aa_\beta(\xi)$.\hfill$\Box$
\end{proof}

\subsection{Integration over non-tangential approach regions}
For $F$ a closed subset of $\S^{n-1}$, the tent over $F$
is defined by

$$\rr_{\alpha}(F)=\bigcup_{\xi\in F}\aa_\alpha(\xi).$$
We will need the two following lemmas for integration over tents.
Their proofs are similar to the ones for
integration over tents in $\R^{n+1}_+$ ({\it see} \cite{ST2}).

\begin{lemmenum}[10] For every $0<\alpha<1$, there exists a constant
$C_\alpha$ such that for every positive function $\Phi$,

$$\int_F\left\{\int_{\aa_\alpha(\xi)}\Phi(x)dx\right\}d\sigma(\xi)\leq
C_\alpha\int_{\rr_\alpha(F)}\Phi(x)(1-\abs{x})^{n-1}dx.$$
\end{lemmenum}

As in the $\R_+^{n+1}$ case, the converse of this lemma is more 
complicated : let $F$ be a
closed subset of $\S^{n-1}$ and let $0<\gamma<1$. A point $\xi$ of 
$\S^{n-1}$ is called a {\it $\gamma$-density}
point of $F$ if

$$\frac{\sigma\ent{B(\xi)\cap F}}{\sigma\ent{B(\xi)\cap\S^{n-1}}}
\geq\gamma$$
for every ball $B(\xi)$ centered in $\xi$.

We will denote by $F^*$ the set of $\gamma$-density points of $F$.
The converse of lemma 10 is then :

\begin{lemmenum}[11] Let $0<\alpha<1$. Then there exists
$\gamma$, $0<\gamma<1$ sufficiently near to $1$ such that, for every
closed subset $F$ of $\S^{n-1}$ and every positive function $\Phi$, we 
have

$$\int_{\rr_\alpha(F^*)}\Phi(x)(1-\abs{x})^{n-1}dx\leq C_{\alpha,\gamma}
\int_F\left\{\int_{\aa_\alpha(\xi)}\Phi(x)dx\right\}d\sigma(\xi).$$
\end{lemmenum}

A direct consequence of these two lemmas is the following (see
\cite{CMS}) :

\begin{lemmenum}[12] For $0<p<2$, for $0<\alpha,\beta<1$, there
exists constants $C_1,C_2$ such that for every $\cc^1$ function $u$ on
$\B_n$,

$$C_1\norm{S_\alpha\ent{u}}_{L^p(\S^{n-1})}\leq
\norm{S_\beta\ent{u}}_{L^p(\S^{n-1})}\leq
C_2\norm{S_\alpha\ent{u}}_{L^p(\S^{n-1})}.$$
Similar estimates are valid if we replace $S_\alpha$ by $S_\alpha^N$.
\end{lemmenum}

\subsection{Consequences of the mean value inequalities}

Let $l\in\R$ and $f$ a function defined on
$\B_n$. Define $I_lf$ by

$$I_lf(r\zeta)=\int_0^rf(t\zeta)(1-t)^{l-1}dt,\qquad0<r<1,\
\zeta\in\S^{n-1}.$$

The following lemma is a direct consequence of the mean value
inequalities and its proof follows the main lines of the upper
half-line case ({\it see} \cite{ST1} pages 214--216) or the
$\mm$-harmonic function case in \cite{BBG}.

\begin{lemmenum}[13] For $0<\alpha<\beta<1$, $\gamma>-\frac{n}{2}$,
$l\in\R$ and $d\in\N$, there exists a constant $C$ such that, for every
$\zeta\in\S^{n-1}$,
and for every $\hh$-harmonic function $u$

$$\int_{\aa_\alpha(\zeta)}\ent{I_l(\nabla^dN^ku)}(z)^2(1-\abs{z})^{2\gamma}dz\leq
C\int_{\aa_\beta(\zeta)}\abs{N^ku(z)}^2(1-\abs{z})^{2(l+\gamma-d)}dz.$$
\end{lemmenum}

\begin{remark} If $l$ is a positive integer, then $N^lh=g$
implies

$$\abs{h}\leq C\ent{I_l\abs{g}+\max_{j\leq
l-1,\abs{z}<\eps}\abs{\nabla^jh(z)}}.$$
\end{remark}

This leads to the following lemma ({\it see} \cite{BBG} for the proof in
case of $\mm$-harmonic functions) :

\begin{lemmenum}[14] For $0<\alpha<\beta<1$, $\gamma>-\frac{n}{2}$
and
$d\in\N$, there exists a constant $C$ such that, for every
$\zeta\in\S^{n-1}$, and every $\hh$-harmonic function $u$,

$$\int_{\aa_\alpha(\zeta)}\abs{\nabla^du}(z)^2(1-\abs{z})^{2\gamma}dz
\leq C\int_{\aa_\beta(\zeta)}\!\abs{N^ku(z)}^2
(1-\abs{z})^{2(k+\gamma-d)}dz
+C\sup_{\abs{z}<\eps}\abs{\nabla^{k-1}u(z)}^2.$$
\end{lemmenum}

The last lemma we will need is also similar to the $\R_+^{n+1}$ case
(\cite{ST1}, page  207) and results directly from the mean value
inequality :

\begin{lemmenum}[15] Let $0<\alpha<\beta<1$, There exists a constant
$C$ 
such that for every $\xi\in\S^{n-1}$ and for every $\hh$-harmonic
function $u$
\begin{enumerate}
\item if $\abs{u}\leq 1$ on $\aa_\beta(\xi)$ then 
$\abs{(1-\abs{x}^2)\nabla u}\leq C$ on $\aa_\alpha(\xi)$,

\item if $S_\beta\ent{u}(\xi)\leq 1$ then $\abs{(1-\abs{x}^2)\nabla
u}\leq C$
on $\aa_\alpha(\xi)$.
\end{enumerate}
\end{lemmenum}

\section{Characterization of $\hh^p$ by maximal
functions, area integrals and Littlewood-Paley $g$-functions}
\label{caracterisation}

{\sl In this section we extend the theory of Fefferman-Stein \cite{FS1}
to the $\hh^p$ spaces.}

\subsection{Maximal Characterization of $\hh^p$}

From the definition of $\mm\ent{u}$ and $\mm_\alpha\ent{u}$, it is
obvious that $\mm\ent{u}\leq\mm_\alpha\ent{u}$, in particular, if 
$\mm_\alpha\ent{u}\in L^p(\S^{n-1})$ then $\mm\ent{u}\in L^p(\S^{n-1})$.
The next proposition claims that the converse is true for
$\hh$-harmonic functions as well as for their normal derivatives.

\begin{proposition}[16] For $0<\alpha<1$, $0<p<+\infty$, for
every integer $k\geq0$ and for every $\hh$-harmonic function $u$, the
following are equivalent :

\begin{enumerate}
\item $\mm\ent{N^k u}\in L^p(\S^{n-1})$,
\item $\mm_\alpha\ent{N^k u}\in L^p(\S^{n-1})$.
\end{enumerate}
Moreover, there exists $C=C_{\alpha,p}$ such that for every 
$\hh$-harmonic function $u$,

$$\norm{\mm\ent{N^ku}}_p\leq\norm{\mm_\alpha\ent{N^ku}}_p\leq 
C\norm{\mm\ent{N^ku}}_p.$$
\end{proposition}

\begin{proof}[Proof] According to lemma 7, for
$a\in\aa_\alpha(\zeta)$

$$\abs{N^ku(a)}^{\frac{p}{2}}\leq
C(1-\abs{a})^{-n}\int_{B\bigl(a,2(1-\abs{a})\eps\bigr)}\abs{N^k
u(\omega)}^{\frac{p}{2}}d\omega.$$
Integrating in polar coordinates $\omega=r\eta$, we see that
$\eta\in B\bigl(\zeta,c(1-\abs{a})\bigr)$ and bounding 
$\abs{N^ku(\omega)}$ by $\mm\ent{N^ku}(\zeta)$ we get

\begin{align}
\abs{N^ku(a)}^{\frac{p}{2}}\leq&C(1-\abs{a})^{-n}
\int_{B\bigl(\zeta,c(1-\abs{a})\bigr)\cap\S^{n-1}}
\ent{\mm\ent{N^ku}(\xi)}^{\frac{p}{2}}d\sigma(\xi)
\times\int_{\abs{a}-2(1-\abs{a})\eps}^{\abs{a}+2(1-\abs{a})\eps}r^{n-1}dr\notag\\
\leq&C(1-\abs{a})^{-n+1}
\int_{B\bigl(\zeta,c(1-\abs{a})\bigr)\cap\S^{n-1}}
\ent{\mm\ent{N^ku}(\xi)}^{\frac{p}{2}}d\sigma(\xi).\notag
\end{align}
But
$\sigma\ent{B\bigl(\zeta,c(1-\abs{a})\bigr)\cap\S^{n-1}}\sim(1-\abs{a})^{n-1}$
therefore

$$\mm_\alpha\ent{N^ku(\zeta)}^{\frac{p}{2}}\leq
C\mm_{HL}\ent{\mm\ent{N^ku}^{\frac{p}{2}}}(\zeta)$$
where $\mm_{HL}$ is Hardy-Littlewood's maximal function on
$\S^{n-1}$. We just have to use the fact that $\mm_{HL}$ is bounded
$L^2(\S^{n-1})\mapsto L^2(\S^{n-1})$ to complete the
proof.\hfill$\Box$
\end{proof}

\begin{remarquenum}[1] This proposition, whose proof is directly 
inspired from the $\R^{n+1}_+$ case in \cite{FS1} depends only 
on the mean value inequalities (lemma 7). Thus, it remains 
true if we replace $N^k$ by $\nabla^k$ or by $\lll_{i,j}^k$ (thus also
by $(-\Delta_\sigma)^{k/2}$)  as long as we replace lemma 7 by 
proposition 8 or by inequality (3.5).
\end{remarquenum}

\begin{remarquenum}[2] For $k=0$ this is equivalence
$(1)\Leftrightarrow(2)$ of theorem A.
\end{remarquenum}

\subsection{$\hh_\delta$-harmonic functions}
\label{difficulte}

To take advantage of inequalities on harmonic functions on $\R_+^{n+1}$,
one is often led to introduce the function $u_\eps(x,t)=u(x,t+\eps)$
which is still harmonic if $u$ is, and which is smooth up to the
boundary. One then hopes to get estimates that are independent of 
$\eps$ and then let $\eps$ go to $0$.

In the case of $\hh$-harmonic functions, we would like to operate in the
same
way. Unfortunately, if $u$ is $\hh$-harmonic, the function
$u_\eps(x)=u\bigl((1-\eps)x\bigr)$ may not be $\hh$-harmonic. We are
thus led
to introduce the following notion.

\begin{figure}[ht]
\begin{center}
\includegraphics{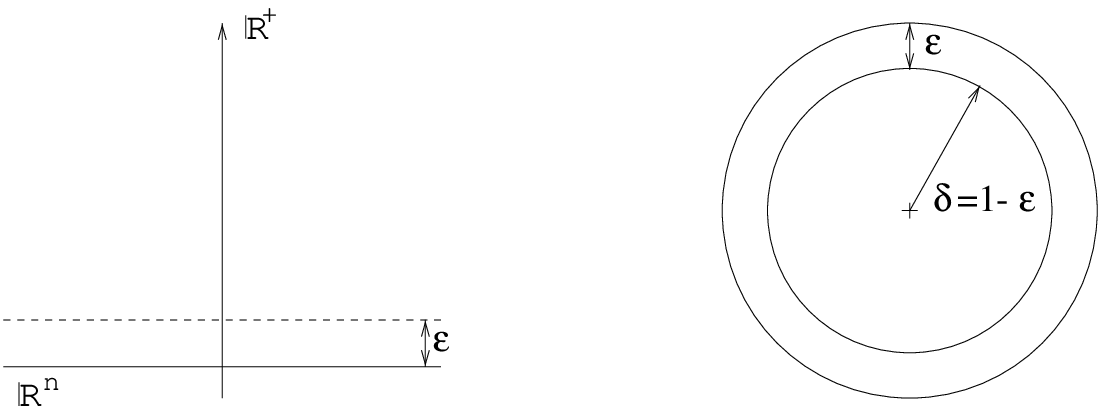}
\end{center}
\caption{Function $u_\eps$}
\end{figure}

\begin{definitionfr} Let $0<\delta<1$ and let $D_\delta$ be the
operator defined by

$$D_\delta=(1-\delta^2r^2)^2\Delta+2(n-2)\delta^2(1-\delta^2r^2)
\sum_{i=1}^nx_i\frac{\partial}{\partial x_i}.$$
We will say that $u$ a smooth function on $\B_n$ is
$\hh_\delta$-harmonic
if $D_\delta u=0$.
\end{definitionfr}

An easy computation shows that if $u$ is $\hh$-harmonic, then the
function
$v$ defined by $v(x)=u(\delta x)$ is $\hh_\delta$-harmonic, {\it i.e.}
$D_\delta u=0$ or also $L_\delta u=0$ with

$$L_\delta=(1-\delta^2r^2)\Delta+2(n-2)\delta^2
\sum_{i=1}^nx_i\frac{\partial}{\partial x_i}.$$
To this laplacian, one can associate the $\hh_\delta$-Poisson kernel
given by its spherical harmonics expansion

$$\P_{h,\delta}(r\zeta,\xi)=\sum_{l=0}^\infty
f_l(\delta^2r^2)r^l\Z^l_\zeta(\xi),$$
to which we can associate $\hh_\delta$-Poisson integrals.

Note also that we obtain a family of operators
$D_\delta$ such that $D_0=\Delta$, the Euclidean laplacian and $D_1=D$ 
the hyperbolic laplacian. Similarely, notice that
$\P_{h,0}=\P_e$ and $\P_{h,1}=\P_h$.

Green's formula for this Laplacian is

$$\int_\Omega\bigl(u(x)D_\delta v(x)
-v(x)D_\delta u(x)\bigr)\frac{dx}{(1-\delta^2\abs{x}^2)^{n}}
=\int_{\partial \Omega}\ent{u\frac{\partial v}{\partial
\vec{n}}-v\frac{\partial u}{\partial\vec{n}}}
\frac{d\sigma}{(1-\delta^2r^2)^{n-2}}.$$

One can check that proofs of chapter \ref{inmoy}\ can be reproduced for 
$\hh_\delta$-harmonic functions, in particular, the mean value
inequality (lemma 7) remains valid with constants independent from
$\delta$. More precisely, we obtain :

\begin{lemmenum}[17] For every
$\eps<\frac{1}{6}$, $k,d\in\N$, $0<p<+\infty$, there exists a constant
$c$
such that, for every $\frac{1}{2}<\delta\leq1$, every
$\hh_\delta$-harmonic
function $u$ and every $a\in\B_n$,

$$\abs{\nabla^dN^ku(a)}\leq c(1-\abs{a})^{-d-\frac{n}{p}}\left(
\int_{B(a,6(1-\abs{a}^2)\eps)}\abs{N^ku(x)}^pdx\right)^{\frac{1}{p}}.$$
\end{lemmenum}

We will need the following estimates :

\begin{proposition}[18] There exists a constant $C$ such that for
every $0<\delta<1$, every $x\in\B_n$ and every $\xi\in\S^{n-1}$,

$$0\leq\P_{h,\delta}(x,\xi)\leq C\P_e(x,\xi).$$
In particular, for every $0<\alpha<1$, there exists a constant $C$ such
that if $f\in L^2(\S^{n-1})$ and $u=\P_{h,\delta}\ent{f}$
then 

$$\norm{\mm_\alpha\ent{u}}_{L^2(\S^{n-1})}\leq
C\norm{f}_{L^2(\S^{n-1})}.$$

Conversely, there exists $\eta>0$ and a constant $C>0$ such that for all
$0<\delta<1$, for all $\xi\in\S^{n-1}$ and all $x\in\aa_\eta(\xi)$,

$$\P_{h,\delta}(x,\xi)\geq\frac{C}{(1-\abs{x}^2)^{n-1}}.$$
\end{proposition}

\begin{proof}[Proof] $\hh_\delta$-harmonic functions satisfy
the maximum principle, so the Poisson kernel $\P_{h,\delta}$
is positive.

Fix $\xi_0\in\S^{n-1}$ and let $u(x)=\P_{h,\delta}(x,\xi_0)$.
With the mean value inequality (lemma 17)

$$0\leq u(x)\leq\frac{C}{(1-\abs{x}^2)^n}
\int_{B(x,1-\abs{x}^2)}u(y)dy
\leq\frac{C}{(1-\abs{x}^2)^n}
\int_{\abs{x}-(1-\abs{x}^2)}^{\abs{x}+(1-\abs{x}^2)}\int_{\S^{n-1}}
u(r\zeta)d\sigma(\zeta)r^{n-1}dr$$
and as $\int_{\S^{n-1}}u(r\zeta)d\sigma(\zeta)=1$, we get

$$\P_{h,\delta}(x,\xi_0)\leq\frac{C}{(1-\abs{x}^2)^{n-1}},$$
with $C$ independent from $\xi_0$ and from $\delta$. But
$\P_e(x,\xi_0)\simeq\frac{C_\alpha}{(1-\abs{x}^2)^{n-1}}$ in
$\aa_\alpha(\xi_0)$ thus $\P_h\leq C\P_e$ in the interior of 
$\aa_\alpha(\xi_0)$.

On the other hand, $N\P_e(r\zeta,\xi_0)$ has same sign as 

$$-2r\ent{(1-r)^2+r(1-<\zeta,\xi_0>)}
-n(1-r^2)\bigl(r-1+(1-<\zeta,\xi_0>)\bigr)$$
so it is negative in $\B_n\setminus\aa_\alpha(\xi_0)$ for $\alpha$ big
enough. This leads to $D_\delta\P_e<0$ on 
$\B_n\setminus\aa_\alpha(\xi_0)$ and
$C\P_e(x,\xi_0)-\P_h(x,\xi_0)>0$ on the boundary of
$\B_n\setminus\aa_\alpha(\xi_0)$ (with $C$ independent from $\delta$ 
and from $\xi_0$), thus, by the maximum principle $\P_{h,\delta}\leq
C\P_e$ on
$\B_n\setminus\aa_\alpha(\xi_0)$.\hfill$\diamond$

For the other inequality, first notice that

$$\P_{h,\delta}(r\zeta,\xi)=\sum_{l\geq0}
f_l(r^2)r^l\Z^l_\zeta(\xi),$$
and as $\Z^l_\xi(\xi)=1$, it turns out that

$$\P_{h,\delta}(r\xi,\xi)\geq\P_e(r\xi,\xi)\geq\frac{C_1}{(1-r^2)^{n-1}}.$$

But $\P_{h,\delta}$ is $\hh_\delta$-harmonic and therefore satisfies
mean value inequalities (lemma 17), {\it i.e.}

\begin{align}
\abs{\nabla\P_{h,\delta}(x,\xi)}\leq&\frac{C}{(1-\abs{x}^2)^{n+1}}
\int_{B\bigl(x,6(1-\abs{x}^2)\eps\bigr)}\abs{\P_{h,\delta}(y,\xi)}dy\notag\\
\leq&\frac{C}{(1-\abs{x}^2)^{n+1}}
\int_{B\bigl(x,6(1-\abs{x}^2)\eps\bigr)}\abs{\P_e(y,\xi)}dy\notag\\
\leq&\frac{C_2}{(1-\abs{x}^2)^n}.\notag
\end{align}
Thus, by the fundamental theorem of calculus,

\begin{align}
\P_{h,\delta}(x,\xi)\geq&\P_h(\abs{x}\xi,\xi)-d(x,\abs{x}\xi)\sup_{\ent{x,\abs{x}\xi}}
\abs{\nabla\P_{h,\delta}(x,\xi)}\notag\\
\geq&\frac{C_1}{(1-\abs{x}^2)^{n-1}}-
C_2\frac{d(x,\abs{x}\xi)}{(1-\abs{x}^2)^n}.\notag
\end{align}
Then, if $\eta$ is small enough to have $C_2d(x,\abs{x}\xi)\leq
C_1(1-\abs{x}^2)$
in $\aa_\eta(\xi)$, then in $\aa_\eta(\xi)$,

$$\P_{h,\delta}(x,\xi)\geq\frac{C}{(1-\abs{x}^2)^{n-1}}.\eqno\Box$$
\end{proof}

\subsection{Characterization by area integral}

In this chapter we characterize $\hh^p$ in terms of area integrals. The
proof is inspired by \cite{FS1} but needs an adaptation to the fact
that $\hh_\delta$-harmonic functions are not $\hh$-harmonic. More
precisely, we will prove the following part of theorem A :

\begin{theorem}[19] For $0<p<2$ and $u$ $\hh$-harmonic, the
following are equivalent :
\begin{enumerate}
\item $\mm_\alpha\ent{u}\in L^p$ for some $\alpha$, $0<\alpha<1$,
\item $S_\alpha\ent{u}\in L^p$ for some $\alpha$, $0<\alpha<1$,
\item $S_\alpha^N\ent{u}\in L^p$ for some $\alpha$, $0<\alpha<1$.
\end{enumerate}
\end{theorem}

\begin{proof}[Proof] Let us show that
$\norm{S_\beta\ent{u}}_p\leq C\norm{\mm_\alpha\ent{u}}_p$. 
According to proposition 16, we may assume
$\alpha<\beta$. Assume first that $u$ is the Poisson integral of an
$L^2$ function.

\begin{notation} For a  measurable function $f:\S^{n-1}\mapsto\R$, we will write 

$$\lambda_f(x)=\sigma\ent{\{\xi\in\S^{n-1}\ :\ \abs{f(\xi)}>x\}}.$$
\end{notation}

Fix momentarily $\mu>0$. Let 
$E=\{\xi\in\S^{n-1}\ :\ \mm_\alpha\ent{u}\leq\mu\}$ and
$B=\S^{n-1}\setminus E$ so that
$\lambda_{\mm_\alpha\ent{u}}(\mu)=\sigma(B)$.
Let $\rr=\rr_\beta(E)=\bigcup_{\xi\in E}\aa_\beta(\xi)$. There exists
an increasing sequence of domains $\rr_\eps$ with $\cc^1$ smooth
boundary
such that $\rr_\eps\rightarrow\rr$.

\begin{figure}[ht]
\begin{center}
\includegraphics{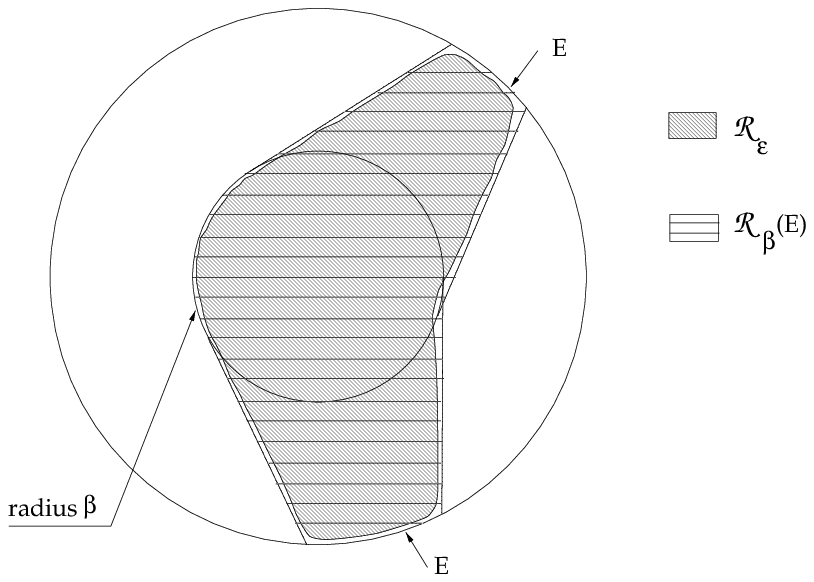}
\end{center}
\caption{$\rr_\beta(F)$ et $\rr_\eps$}
\end{figure}

\noindent Further, by definition of $E$, we have $\abs{u(x)}\leq\mu$ in
$\bigcup_{\xi\in E}\aa_\alpha(\xi)$ thus, with lemma 15, 
$\abs{(1-\abs{x}^2)\nabla u}\leq C\mu$ in $\rr$. Then

\begin{align}
\int_ES_\beta\ent{u}^2(\xi)d\sigma(\xi)=&\int_E\ent{\int_{\aa_\beta(\xi)}
\abs{\nabla u(x)}^2(1-\abs{x}^2)^{-n+2}dx}d\sigma(\xi)\notag\\
\leq&C\int_\rr\abs{\nabla u(x)}^2(1-\abs{x}^2)dx\notag
\end{align}
according to lemma 10. But $L\abs{u}^2=2(1-\abs{x}^2)\abs{\nabla
u(x)}^2$, therefore

\begin{align}
\int_ES_\beta\ent{u}^2(\xi)d\sigma(\xi)\leq&C\int_\rr
L\abs{u}^2(1-\abs{x}^2)^{n-1}
\frac{dx}{(1-\abs{x}^2)^{n-1}}\notag\\
\leq&C\overline{\lim_{\eps\rightarrow0}}\int_{\rr_\eps}L\abs{u}^2(1-\abs{x}^2)^{n-1}
\frac{dx}{(1-\abs{x}^2)^{n-1}}\tag{4.1}
\end{align}
Since $L(1-\abs{x}^2)^{n-1}=-2(n-1)(1-\abs{x}^2)^{n-1}$,

\begin{align}
\int_{\rr_\eps}L\abs{u}^2(1-\abs{x}^2)^{n-1}
\frac{dx}{(1-\abs{x}^2)^{n-1}}
=&\int_{\rr_\eps}\ent{(1-\abs{x}^2)^{n-1}L\abs{u}^2
-\abs{u}^2L(1-\abs{x}^2)^{n-1}}\frac{dx}{(1-\abs{x}^2)^{n-1}}\notag\\
&-2(n-1)\int_{\rr_\eps}\abs{u}^2dx\notag\\
\leq&\int_{\rr_\eps}\ent{(1-\abs{x}^2)^{n-1}L\abs{u}^2
-\abs{u}^2L(1-\abs{x}^2)^{n-1}}\frac{dx}{(1-\abs{x}^2)^{n-1}}\notag
\end{align}
and Green's formula leads to

$$\int_ES_\beta\ent{u}^2(\xi)d\sigma(\xi)\leq
C\overline{\lim_{\eps\rightarrow0}}\int_{\partial \rr_\eps}\ent{
(1-\abs{x}^2)^{n-1}\frac{\partial\abs{u}^2}{\partial n}-\abs{u}^2
\frac{\partial(1-\abs{x}^2)^{n-1}}{\partial n}}
\frac{d\sigma_{\rr_\eps}(x)}{(1-\abs{x}^2)^{n-2}}.$$

Now cut $\partial\rr_\eps$ into two parts
$\partial\rr_\eps^E$ and
$\partial\rr_\eps^B$ where

$$\partial\rr_\eps^E=\{r\xi\in\partial\rr_\eps\ :\ \xi\in E\}\quad
\mathrm{and}\quad
\partial\rr_\eps^B=\{r\xi\in\partial\rr_\eps\ :\ \xi\in B\}.$$

--- On $\partial\rr_\eps^E$, we have $d\sigma_{\rr_\eps}\sim d\sigma$.
Further $\sup_{r>0}\abs{u(r\xi)}$ and 
$\sup_{r>0}(1-r^2)\abs{\nabla u(r\xi)}$ are in $L^2(\S^{n-1})$. Finaly,
as $u$ is the Poisson integral of an $L^2$ function,
$\lim_{r\rightarrow1}(1-r^2)\abs{\nabla u(r\xi)}=0$ almost
everywhere thus, by Lebesgues' lemma,

$$\int_{\partial\rr_\eps^E}(1-\abs{x}^2)\frac{\partial\abs{u}^2}{\partial
n}
d\sigma_{\rr_\eps}(x)\rightarrow0$$
when $\eps\rightarrow0$.

--- We have already seen that $\abs{(1-\abs{x}^2)\nabla u(x)}\leq C\mu$
and that $\abs{u(x)}\leq\mu$ in $\rr$, thus

$$\int_{\partial\rr_\eps^B}(1-\abs{x}^2)\frac{\partial\abs{u}^2}{\partial
n}
d\sigma_{\rr_\eps}(x)\leq
C\mu^2\int_{\partial\rr_\eps^B}d\sigma_{\rr_\eps}
\leq C\mu^2\sigma(B)=C\mu^2\lambda_{\mm_\alpha\ent{u}}(\mu).$$

--- On the other hand

$$\abs{\int_{\partial\rr_\eps^E}\abs{u}^2
\frac{\partial(1-\abs{x}^2)^{n-1}}{\partial n}
\frac{d\sigma_{\rr_\eps}(x)}{(1-\abs{x}^2)^{n-2}}}
\leq C\int_E\mm_\alpha\ent{u}^2(\xi)d\sigma(\xi)
\leq C\int_0^\mu t\lambda_{\mm_\alpha\ent{u}}(t)dt$$
since $\mm_\alpha\ent{u}\leq\mu$ on $E$.

--- Finally

$$\abs{\int_{\partial\rr_\eps^B}\abs{u}^2
\frac{\partial(1-\abs{x}^2)^{n-1}}{\partial n}
\frac{d\sigma_{\rr_\eps}(x)}{(1-\abs{x}^2)^{n-2}}}
\leq C\mu^2\int_Bd\sigma=C\mu^2\sigma(B)
= C\mu^2\lambda_{\mm_\alpha\ent{u}}(\mu).$$

But then

\begin{align}
\lambda_{S_\beta\ent{u}}(\mu)=&\sigma\ent{\{\xi:\abs{S_\beta(u)(\xi)}>\mu\}}\notag\\
=&\sigma\ent{\{\xi:\abs{S_\beta(u)(\xi)}>\mu,\text{ and }
\mm_\alpha\ent{u}(\xi)\leq\mu\}}
+\sigma\ent{\{\xi:\abs{S_\beta(u)(\xi)}>\mu,\text{ and }
\mm_\alpha\ent{u}(\xi)>\mu\}}\notag\\
\leq&\frac{1}{\mu^2}\int_{\mm_\alpha\ent{u}\leq\mu}S_\beta\ent{u}^2(\xi)
d\sigma(\xi)+\sigma\ent{\{\xi\ :\ \mm_\alpha\ent{u}(\xi)>\mu\}}\notag\\
\leq&\frac{1}{\mu^2}\int_ES_\beta\ent{u}^2(\xi)d\sigma(\xi)
+\lambda_{\mm_\alpha\ent{u}}(\mu)\notag
\end{align}
and taking into account the previous estimates, we get

$$\lambda_{S_\beta\ent{u}}(\mu)\leq
C\ent{\lambda_{\mm_\alpha\ent{u}}(\mu)
+\frac{1}{\mu^2}\int_0^\mu t\lambda_{\mm_\alpha\ent{u}}(t)dt}.
\leqno(4.2)$$
After integrating, we get

\begin{align}
\norm{S_\beta\ent{u}}_p^p=p\int_0^\infty\mu^{p-1}\lambda_{S_\beta\ent{u}}(\mu)
d\mu
\leq&C\int_0^\infty\mu^{p-1}\lambda_{\mm_\alpha\ent{u}}(\mu)
d\mu+C\int_0^\infty\mu^{p-2}\int_0^\mu
t\lambda_{\mm_\alpha\ent{u}}(t)dt
d\mu\notag\\
\leq&C\norm{\mm_\alpha(u)}^p_p+C\int_0^\infty
t\lambda_{\mm_\alpha\ent{u}}(t)
\int_t^\infty\mu^{p-3}d\mu dt\notag\\
\leq&C\norm{\mm_\alpha(u)}^p_p\notag
\end{align}
since $0<p<2$.

We have shown that $\norm{S_\beta\ent{u}}_{L^p}\leq
C\norm{u}_{\hh^p}$ for every $u\in\hh^p\cap
\P_h\ent{L^2(\S^{n-1})}$. But, with help of the atomic decomposition of
$\hh^p$ ({\it see} \cite{Ja3}),
$\P_h\ent{L^2(\S^{n-1})}$ is dense in $\hh^p$, we deduce the inequality
for
every $u\in\hh^p$.\hfill$\diamond$

Let us now show the implication ``$(2)\Rightarrow(1)$'' for $0<p<2$.
More precisely, we will show that there exists a constant $C$ such that
for every $\hh$-harmonic function $u$, $\norm{\mm_\alpha\ent{u}}_p
\leq C\norm{S_\beta\ent{u}}_p$.

With proposition 16, up to a change of the constant $C$, we may
assume that $\alpha<\beta$ and $\alpha<\eta$ where $\eta$ is given by
proposition 18 to be such
that $\P_{h,\delta}(x,\xi)\geq\frac{C}{(1-\abs{x}^2)^{n-1}}$ on
$\aa_\eta(\xi)$, an estimate we will use at the end of the proof of the 
theorem (see the proof of the claim).
Let $\frac{1}{2}<\delta<1$. Let $u$ be an $\hh$-harmonic function
satisfying $S_\beta\ent{u}\in L^p(\S^{n-1})$, and let
$u_\delta(x)=u(\delta x)$, 
in particular $u_\delta$ is a $\cc^\infty$ function on $\overline{\B_n}$.

We will show that
$\norm{\mm_\alpha\ent{u_\delta}}_p\leq C\norm{S_\beta\ent{u_\delta}}_p$
with $C$ a constant independent from $\delta$. The result follows by
letting $\delta\rightarrow1$.

For $\mu>0$, let $E=\{\xi\in\S^{n-1}\ :\
S_\beta\ent{u_\delta}(\xi)\leq\mu\}$
and let $B=\S^{n-1}\setminus E$, therefore 
$\lambda_{S_\beta\ent{u_\delta}}(\mu)=\sigma(B)$. Let $E_0$ be the set
of
$\gamma$-density points of $E$ where $\gamma$ is chosen so as to be
able to apply lemma 11. Let $B_0=\S^{n-1}\setminus E_0$.
Note that by Hardy-Littlewood's maximal theorem,
$\sigma(B_0)\leq C\sigma(B)\leq C\lambda_{S_\beta\ent{u}}(\mu)$.

Let $\rr=\bigcup_{\xi\in E_0}\aa_\alpha(\xi)$ and let $\rr_\eps$ be a
sequence
of domains with $\cc^1$ boundary approximating $\rr$ and such that
$dist(\rr_\eps,\S^{n-1})\geq\eps$. We have

\begin{align}
\int_{E_0}S_\beta\ent{u_\delta}(\xi)^2d\sigma(\xi)=&\int_{E_0}\ent{
\int_{\aa_\beta(\xi)}\abs{\nabla u_\delta(x)}^2(1-\abs{x}^2)^{-n+2}dx}
d\sigma(\xi)\notag\\
\geq&C\int_\rr(1-\abs{x}^2)\abs{\nabla u_\delta(x)}^2dx\notag\\
\geq&C\int_{\rr_\eps}(1-\abs{x}^2)\abs{\nabla u_\delta(x)}^2dx\notag
\end{align}
according to lemma 11.

Write $\displaystyle I_\delta=\int_{\rr_\eps}(1-\abs{x}^2)
\abs{\nabla u_\delta(x)}^2dx$. As $u_\delta$ is $\hh_\delta$-harmonic,
$L_\delta\abs{u_\delta}^2(x)=2(1-\delta^2\abs{x}^2)\abs{\nabla
u_\delta(x)}^2$. Write $v(x)=(1-\abs{x}^2)(1-\delta^2\abs{x}^2)^{n-2}$,
so that

\begin{align}
I_\delta=&\!\int_{\rr_\eps}
(1-\abs{x}^2)(1-\delta^2\abs{x}^2)^{n-2}L_\delta\abs{u_\delta}^2(x)\!
-\!\abs{u_\delta}^2(x)L_\delta v(x)
\!\frac{dx}{(1-\delta^2\abs{x}^2)^{n-1}}\notag\\
&+\int_{\rr_\eps}\abs{u_\delta(x)}^2\frac{L_\delta
v(x)}{(1-\delta^2\abs{x}^2)^{n-1}}dx.\notag
\end{align}
Call $\tilde\aa_\beta(\xi)=\aa_\beta(\xi)\cap\rr_\eps$, as

\begin{align}
L_\delta v(x)=&-2n(1+(n-2)\delta^2-\delta^2(n-1)\abs{x}^2)
(1-\delta^2\abs{x}^2)^{n-2}
+4\abs{x}^2(n-2)(1-\delta^2)(1-\delta^2\abs{x}^2)^{n-3}\notag\\
\geq&-2n\bigl(1+(n-2)\delta^2\bigr)
\ent{1-\frac{\delta^2(n-1)}{1+(n-2)\delta^2}\abs{x}^2}\geq-C\notag
\end{align}
we get

\begin{align}
\int_{\rr_\eps}\abs{u_\delta}^2
\frac{L_\delta v(x)}{(1-\delta^2\abs{x}^2)^{n-1}}dx\geq&
-C\int_{E_0}\int_{\tilde\aa_\beta(\xi)}
\abs{u_\delta(x)}^2\frac{1}{(1-\delta^2\abs{x}^2)^{n-1}}dx\notag\\
\geq&-C\int_{E_0}\int_{\tilde\aa_\beta(\xi)}\abs{I_1\nabla
u_\delta(x)}^2
\frac{dx}{(1-\delta^2\abs{x}^2)^{n-1}}\notag\\
\geq&-C\int_{E_0}\int_{\aa_{\beta'}(\xi)}\abs{\nabla u_\delta(x)}^2
(1-\abs{x}^2)^2\frac{dx}{(1-\abs{x}^2)^{n-1}}\notag\\
\geq&-C\int_{E_0}S_{\beta'}\ent{u_\delta}(\xi)^2d\sigma(\xi)\notag
\end{align}
with lemma 13, for some $\beta'>\beta$. Finally

$$\int_{E_0}S_{\beta'}\ent{u_\delta}(\xi)^2d\sigma(\xi)\geq CJ_\delta$$
where

$$J_\delta=\int_{\rr_\eps}
(1-\abs{x}^2)(1-\delta^2\abs{x}^2)^{n-2}L_\delta\abs{u_\delta}^2
-\abs{u_\delta}^2L_\delta v(x)\frac{dx}{(1-\delta^2\abs{x}^2)^{n-1}}.$$
Green's formula then leads to

\begin{align}
J_\delta=&\int_{\partial\rr_\eps}\!(1-\abs{x}^2)
\frac{\partial\abs{u_\delta}^2}{\partial\vec{n}}(x)
-\abs{u_\delta}^2(x)(1-\abs{x}^2)^{2-n}\frac{\partial
v(x)}{\partial\vec{n}}
d\sigma_{\rr_\eps}\notag\\
\geq&C_1\int_{\partial\rr_\eps}\!\abs{u_\delta}^2d\sigma_{\rr_\eps}
\!-C_2\!\int_{\partial\rr_\eps}\!(1-\abs{x}^2)\abs{u_\delta}
\abs{\nabla u_\delta}d\sigma_{\rr_\eps}
\tag{4.3}
\end{align}
since $-\frac{\partial v}{\partial\vec{n}}\geq C_1(1-r^2)^{n-2}$ and 
$C_1,C_2$ are independent from $\eps$ and from $\delta$.

Let $K_\eps=\left(\int_{\partial\rr_\eps}\abs{u_\delta(x)}^2
d\sigma_{\rr_\eps}(x)\right)^{\frac{1}{2}}$ which is finite since
$u_\delta$ is $\cc^\infty$.

Then again, cut $\partial\rr_\eps$ into two parts,
$\partial\rr_\eps=\partial\rr_\eps^E\cup\partial\rr_\eps^B$ with

$$\partial\rr_\eps^E=\{r\xi\in\partial\rr_\eps\ :\ \xi\in E_0\}\quad
\mathrm{and}\quad
\partial\rr_\eps^B=\{r\xi\in\partial\rr_\eps\ :\ \xi\in B_0\}.$$

With lemma 15, $\abs{(1-\abs{x}^2)\nabla u_\delta(x)}\leq C\mu$ 
in $\rr$ since $S_\beta\ent{u_\delta}(\xi)\leq\mu$ in $E_0$ ($C$
independent from $\delta$). The Cauchy-Schwartz inequality then gives

\begin{align}
\int_{\partial\rr_\eps^B}(1-\abs{x}^2)\abs{u_\delta(x)}
\abs{\nabla u_\delta(x)}d\sigma_{\rr_\eps}(x)\leq&
K_\eps\mu\sigma_{\rr_\eps}(\partial\rr_\eps^B)^{\frac{1}{2}}\leq 
CK_\eps\mu\sigma(B)^{\frac{1}{2}}\notag\\
\leq&CK_\eps\bigl(\mu^2\lambda_{S_\beta\ent{u_\delta}}(\mu)
\bigr)^{\frac{1}{2}}.\notag
\end{align}
As $u_\delta$ is $\cc^\infty$,

$$M_\eps=\int_{\partial\rr_\eps^E}(1-\abs{x}^2)\abs{u_\delta}
\abs{\nabla u_\delta}d\sigma_{\rr_\eps}\rightarrow 0$$
when $\eps\rightarrow0$, thus (4.3) and $J_\delta\leq
C\int_{E_0}S_\beta\ent{u_\delta}(\xi)^2d\sigma(\xi)$ imply

$$K_\eps^2\leq C\int_{E_0}S_\beta\ent{u_\delta}(\xi)^2d\sigma(\xi)+
CK_\eps\bigl(\mu^2\lambda_{S_\beta\ent{u_\delta}}(\mu)\bigr)^{\frac{1}{2}}
+CM_\eps$$
thus, if $\eps$ is small enough,

$$K_\eps^2\leq C\ent{\int_{E_0}S_\beta\ent{u_\delta}(\xi)^2d\sigma(\xi)+
\mu^2\lambda_{S_\beta\ent{u_\delta}}(\mu)}.
\leqno(4.4)$$
Next, define
$f_\eps(\xi)=\abs{u_\delta\bigl(r_\eps(\xi)\xi\bigr)}
+\mu\chi_{B_0}(\xi)$ where $r_\eps(\xi)\xi$ is a 
parameterization of $\partial\rr_\eps$. In vue of (4.4),
$f_\eps$ is an $L^2$ function and

$$\left(\int_{\S^{n-1}}\abs{f_\eps(\xi)}^2d\sigma(\xi)
\right)^{\frac{1}{2}}\leq 2K_\eps+2\mu\sigma(B)^{\frac{1}{2}}.$$
Let $U_\eps(x)$ be its $\hh_\delta$-Poisson integral.

{\bf Claim :} {\sl $\abs{u_\delta(x)}\leq CU_\eps(x)$ in $\rr_\eps$.}

We postpone the proof of this claim to the end of the proof
of the theorem.

Taking a subsequence of $f_\eps$ that converges weakly to a
function $f\in L^2$, it results from (4.4) that

$$\int_{\S^{n-1}}\abs{f(\xi)}^2d\sigma(\xi)\leq 
C\ent{\int_{E_0}S_\beta\ent{u_\delta}(\xi)^2d\sigma(\xi)+
\mu^2\lambda_{S_\beta\ent{u_\delta}}(\mu)}.$$
On the other hand, as $\abs{u_\delta(x)}\leq U_\eps(x)$ in $\rr_{\eps}$,
going to the limit, $\abs{u_\delta(x)}\leq U(x)$ in $\rr$ where 
$U=\P_{h,\delta}\ent{f}$ is the $\hh_\delta$-Poisson integral of $f$,
thus
for $x\in E_0$, $\mm_\alpha\ent{u_\delta}(x)\leq C
\mm_\alpha\ent{U}(x)$. So

$$\int_{E_0}\mm_\alpha\ent{u_\delta}(\xi)^2d\sigma(\xi)\leq
C\int_{E_0}\mm_\alpha\ent{U}(\xi)^2d\sigma(\xi)
\leq C\int_{\S^{n-1}}\abs{f(\xi)}^2d\sigma(\xi).$$
It follows that

$$\sigma\{\xi\in E_0:\mm_\alpha\ent{u_\delta}\geq\mu\}\leq 
C\ent{\lambda_{S_\beta}\ent{u_\delta}(\mu)+\frac{1}{\mu^2}\int_0^\mu
t\lambda_{S_\beta\ent{u_\delta}}(t)dt}.$$
and as $\sigma(\S^{n-1}\setminus E_0)=\sigma(B)\leq
 C\lambda_{S_\beta\ent{u_\delta}}(\mu)$, we get

$$\lambda_{\mm_\alpha\ent{u_\delta}}(\mu)\leq 
C\ent{\lambda_{S_\beta\ent{u_\delta}}(\mu)+
\frac{1}{\mu^2}\int_0^\mu t\lambda_{S_\beta\ent{u_\delta}}(t)dt}$$
and an integration similar to the one after inequality
(4.2)
implies that there exists a constant $C$ such that

$$\norm{\mm_\alpha\ent{u_\delta}}_p\leq
C\norm{S_\beta\ent{u_\delta}}_p.$$
We conclude by letting $\delta$ go to $1$.\hfill$\diamond$

{\sl Proof of claim :} Let
$f_\eps(x)=\abs{u_\delta(r_\eps(\xi)\xi)}+\mu\chi_{B_0}(\xi)$ et
$U_\eps(x)=\P_{h,\delta}\ent{f_\eps}(x)$.

We want to show that, for $x\in\rr_\eps$, $\abs{u_\delta(x)}\leq
CU_\eps(x)$. By the maximum principle, it is enough to prove this
inequality on
$x\in\partial\rr_\eps=\partial\rr_\eps^E\cup\partial\rr_\eps^B$.

--- On $\partial\rr_\eps^E$, the inequality is true as long as we take $C$
big enough.

--- Recall that on $\rr_\eps$, $(1-\abs{x})\abs{\nabla
u_\delta(x)}\leq C\mu$. Then, if
$x_1,x_2\in\rr_\eps$, 

$$\abs{u_\delta(x_1)-u_\delta(x_2)}\leq\abs{x_1-x_2}\sup_{x\in\ent{x_1,x_2}}\abs{\nabla
u_\delta(x)}.$$
Thus, if $x_1\in\rr_\eps$ and if $x_2\in
B\bigl(x_1,\frac{1}{2}(1-\abs{x_1}^2)\bigr)$,
then

$$\abs{u_\delta(x_1)-u_\delta(x_2)}\leq C\mu.
\leqno(4.5)$$
Fix $x_1\in\partial\rr_\eps^B$ and let $\ss_\eps$ be the portion of
$\partial\rr_\eps^B$ located in the ball
$B\bigl(x_1,\frac{1}{2}(1-\abs{x_1}^2)\bigr)$. By (4.5),

$$\abs{u_\delta(x_1)}\leq\frac{1}{\sigma_\eps(\ss_\eps)}\int_{\ss_\eps}
(\abs{u_\delta(x_2)}+C\mu )d\sigma_\eps(x_2).$$
Since $d\sigma_\eps\simeq d\sigma$ and since
$B\bigl(x_1,\frac{1}{2}(1-\abs{x_1}^2)\bigr)$ is of radius
$\frac{1}{2}(1-\abs{x_1}^2)$, $\sigma_\eps(\ss_\eps)\simeq
a(1-\abs{x_1}^2)^{n-1}$. Therefore, by definition of $f_\eps$,

$$\abs{u_\delta(x_1)}\leq\frac{C}{(1-\abs{x_1}^2)^{n-1}}
\int_{E_0}f_\eps(\xi)d\sigma(\xi).$$

But, the Poisson kernel $\P_{h,\delta}(x_1,\xi)$ is
$\geq\frac{c}{(1-\abs{x_1}^2)^{n-1}}$ in $\aa_\alpha(\xi)$, it follows
that

$$\abs{u_\delta(x_1)}\leq
C\int_{E_0}\P_{h,\delta}(x_1,\xi)f_\eps(\xi)d\sigma(\xi),$$
and thus we have $\abs{u_\delta}\leq C\P_{h,\delta}\ent{f_\eps}$
in $\rr_\eps$.\hfill$\diamond$

The equivalence ``$(2)\Leftrightarrow(3)$'' results immediately from
lemma 14.\hfill$\Box$
\end{proof}

\begin{remarquenum}[1] The theorem is valid for functions $u$ taking
their
values in a Hilbert space instead of $\R$ or $\C$. The key point is
that
equality $L\abs{u}^2=2(1-\abs{x}^2)\abs{\nabla u}^2$ is valid in Hilbert
spaces.
\end{remarquenum}

\begin{remarquenum}[2] For the proof of
``$(1)\Rightarrow(2)$'', we have used density in $\hh^p$ of Poisson
integrals
of $L^2$ functions obtained with the atomic decomposition. We could also
use $\hh_\delta$-harmonic functions.
\end{remarquenum}

%
%

\subsection[Characterization by $g$ function]{Characterization by Lit\-tlewood-Paley's $g$-function}

Due to the mean value inequality for $\hh$-harmonic functions,
one immediatly gets :

\begin{lemmenum}[20] For every $\alpha$ with $0<\alpha<1$ there
exists a
constant $C$ such that for every $\hh$-harmonic function $u$ and every
$\xi\in\S^{n-1}$,

$$g\ent{u}(\xi)\leq CS_\alpha\ent{u}(\xi).$$
\end{lemmenum}

\begin{proof}[Proof] Simply adapt the $\R^{n+1}_+$ case from
\cite{ST1}.\hfill$\Box$.
%
%
%
%
\end{proof}

\begin{theorem}[21] Let $0<p<2$. For every $\hh$-harmonic 
function $u$, the following are equivalent :
\begin{enumerate}
\item $g\ent{u}\in L^p(\S^{n-1})$,

\item $g^N\ent{u}\in L^p(\S^{n-1})$,

\item $S_\alpha\ent{u}\in L^p(\S^{n-1})$ for some $\alpha$, $0<\alpha<1$
(thus for every $\alpha$).
\end{enumerate}
\end{theorem}

\begin{proof}[Proof] Let $\hil$ be the Hilbert space defined by

$$\hil=\left\{\ffi:\ent{0,1}\mapsto\C\ :\
\norm{\ffi}_\hil^2=\int_0^1\abs{\ffi(s)}^2(1-s^2)ds<+\infty\right\}.$$

Let $u$ be an $\hh$-harmonic function such that $g\ent{u}\in
L^p(\S^{n-1})$. For $0<s<1$ define $U(r\zeta)=Nu(rs\zeta)$,
then

\begin{align}
\norm{U(r\zeta)}_{\hil}=\norm{s\mapsto
Nu(rs\zeta)}_\hil=&\int_0^1\abs{Nu(rs\zeta)}^2(1-s^2)ds\notag\\
\leq&\int_0^r\abs{\nabla
u(s\zeta)}^2\left(1-\left(\frac{s}{r}\right)^2\right)\frac{ds}{r}\notag\\
=&\frac{1}{r^3}\int_0^r\abs{\nabla u(s\zeta)}^2(r^2-s^2)ds\notag\\
\leq&Cg\ent{u}(\zeta)^2\notag
\end{align}
so 

$$\mm\ent{U}(\xi)=\sup_{0<r<1}\norm{U(r\zeta)}_\hil\leq
Cg\ent{u}(\zeta)\in L^p(\S^{n-1}).$$
According to remark 1 after the proof of theorem 19,

$$\norm{S_\alpha\ent{U}}_p\leq C\norm{\mm\ent{U}}_p\leq
C'\norm{g\ent{u}}_p.$$
Write $S_\alpha\ent{U}(\zeta)$ with the parameterization $r(\xi)\xi$ of
$\partial\aa_\alpha(\zeta)$ :

$$S_\alpha\ent{U}(\zeta)^2=\int_{\S^{n-1}}\int_0^{r(\xi)}\int_0^1
\abs{\nabla Nu(rs\xi)}^2(1-s^2)sds(1-r^2)^{2-n}r^{n-1}drd\sigma(\xi)$$
and, with the change of variables $t=rs$, we get, changing order of
integration

\begin{align}
S_\alpha\ent{U}(\zeta)^2=&\int_{\S^{n-1}}\int_0^{r(\xi)}\int_t^{r(\xi)}
\abs{\nabla N(t\xi)}^2\left(1-\left(\frac{t}{r}\right)^2\right)
\frac{t}{r}(1-r^2)^{2-n}r^{n-2}drdtd\sigma\notag\\
\geq&\int_{\S^{n-1}}\int_0^{r(\xi)}\abs{\nabla Nu(t\xi)}^2
\int_t^{r(\xi)}(r-t)(1-r)^{2-n}dr t^{n-3}dtd\sigma(\xi)\notag
\end{align}
But, if $1-t>2\bigl(1-r(\xi)\bigr)$

$$\int_t^{r(\xi)}(r-t)(1-r)^{2-n}dr=
\int_{1-r(\xi)}^{1-t}s^{2-n}(1-t-s)ds\geq C(1-t)^{4-n}$$
thus there exists $\beta<\alpha$ such that

\begin{align}
S_\alpha\ent{U}(\zeta)^2\geq&\int_{\aa_\beta(\zeta)}\abs{N^2u(t\xi)}^2
(1-t)^{4-n}t^{n-1}dtd\sigma(\zeta)\notag\\
\geq&C\int_{\aa_{\beta'}(\zeta)}\abs{I_1N^2u(t\xi)}^2
(1-t)^{2-n}t^{n-1}dtd\sigma(\zeta)\notag
\end{align}
with $\beta'<\beta$ according to lemma 13, thus
$S_\alpha\ent{U}(\zeta)^2\geq CS_{\beta'}^N\ent{u}(\zeta)$, which
completes
the proof of $(1)\Leftrightarrow(2)$.

The equivalence $(1)\Leftrightarrow(3)$ results directly from lemma
14.\hfill$\Box$
\end{proof}

\section{Characterization of Hardy-Sobolev spaces}

{\sl In this section, we prove theorems B and C.}

In these theorems, that $\mm$ can be replaced by $\mm_\alpha$ is a direct
consequence of the mean value inequality ({\it see} proposition 16 and 
the remarks following it). We will need the following.

\begin{notation} For an integer $k\geq1$, write $\A_k$ for the set of
indices 

$$\A_k=\{(i,j)\in\N\times\N\ :\ 0\leq i\leq k,\ 0\leq j\leq k+1,\ j\
\mathrm{even},\ 1\leq
i+j\leq k+1\}.$$

\begin{figure}[ht]
\begin{center}
\includegraphics{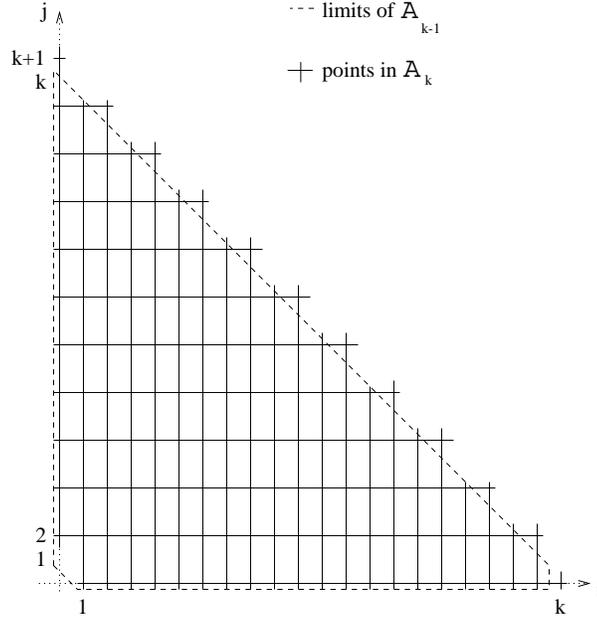}
\end{center}
\caption{The set $A_k$ (here with $k$ even)}
\end{figure}
\end{notation}

\begin{lemmenum}[22] For every $k\geq1$, there exists two families of
polynomials $\bigl(P_j^{(k)}\bigr)_{j=1\ldots\ent{\frac{k}{2}}}$ and 
$\bigl(Q_{i,j}^{(k)}\bigr)_{(i,j)\in\A_k}$ such that for every 
$\hh$-harmonic function $u$ and every $k$,

$$(1-r^2)N^{k+1}u+2(n-1-k)N^ku=
\sum_{j=1}^{\ent{\frac{k}{2}}}P_j^{(k)}(r)\Delta_\sigma^j u
+(1-r^2)\sum_{(i,j)\in\A_k}Q_{i,j}^{(k)}(r)N^i
\Delta_\sigma^{\frac{j}{2}}u.$$
Moreover, the polynomials $P_{\ent{k/2}}^{(k)}$ and 
$Q^{(k)}_{0,\ent{\frac{k+1}{2}}}$ are not zero on the boundary.
\end{lemmenum}

\begin{proof}[Proof] Using the radial-tangential expression of $D$, we
can see that if $Du=0$ then

$$(1-r^2)N^2u+2(n-2)Nu=(1-r^2)\ent{(n-2)Nu-\Delta_\sigma
u}.\leqno(5.1)$$
The lemma is thus verified for $k=1$ with $Q_{1,0}^{(1)}(r)=n-2$,
$Q_{0,2}^{(1)}(r)=-1$.

Applying $N^{k-1}$ to equation (5.1) leads to

\begin{align}
(1-r^2)N^{k+1}u+\bigl(n-1-k\bigr)N^ku=&r^2\sum_{l=1}^{k-1}a_l^{(k)}N^lu
+r^2\sum_{l=0}^{k-2}b_l^{(k)}N^l\Delta_\sigma u\notag\\
&+(1-r^2)\ent{(n-2k)N^ku-N^{k-1}\Delta_\sigma u}\notag
\end{align}
and we conclude with the induction hypothesis.\hfill$\Box$
\end{proof}

The equivalence of $1a$ and $1b$ as well as the equivalence of $2a$ and
$2b$ in theorem C have already been shown. We will now prove the
remaining of this theorem.

\begin{theorem}[23] For $0<\alpha<1$, $0<p<+\infty$, for every 
integer $0\leq k\leq n-2$ and for every $\hh$-harmonic function $u$, the
following are equivalent :

\begin{enumerate}
\item If $k$ is {\rm even}
\begin{enumerate}
\item $\mm_\alpha\ent{N^j u}\in L^p(\S^{n-1})$, for 
$0\leq j\leq k$,
\item $\mm_\alpha\ent{\Delta_\sigma^j u}\in L^p(\S^{n-1})$, for 
$0\leq j\leq \frac{k}{2}$,
\item $\mm_\alpha\ent{\nabla^ju}\in L^p(\S^{n-1})$ for 
$0\leq j\leq k$.
\end{enumerate}
\item If $k$ is {\rm odd}
\begin{enumerate}
\item $\mm_\alpha\ent{N^j u}\in L^p(\S^{n-1})$, for 
$0\leq j\leq k$,
\item $\mm_\alpha\ent{\Delta_\sigma^j u}\in L^p(\S^{n-1})$, for
$0\leq j\leq \frac{k-1}{2}$, and 
$\mm_\alpha\ent{(1-r^2)\Delta_\sigma^{\frac{k+1}{2}}u}\in
L^p(\S^{n-1})$.
\end{enumerate}
\end{enumerate}
\end{theorem}

\begin{remark} If $\mm_\alpha\ent{\Delta_\sigma^{\frac{k}{2}}u}\in
L^p(\S^{n-1})$ then, by the mean value properties, we get that

$$\mm_\alpha\ent{(1-r^2)\Delta_\sigma^{\frac{k+1}{2}}u}\in
L^p(\S^{n-1}).$$
Hence $\hh^p_k\subset\hh^p_{k,N}$.
\end{remark}

\begin{proof}[Proof] The theorem is of course true for $k=0$.
Assume the result holds up to rank $k-1$.

Assume first that $k$ is {\it even}. Implication $(c)\Rightarrow(b)$
is obvious. Let us show $(b)\Rightarrow(a)$.

Let $u$ be an $\hh$-harmonic function that satisfies $(b)$, and let

$$v=2(n-1-k)N^ku+(1-r^2)N^{k+1}u.\leqno(5.2)$$
According to lemma 22,

$$v(r\zeta)=\sum_{j=1}^{\ent{\frac{k}{2}}}P_j^{(k)}(r)\Delta_\sigma^j u
+(1-r^2)\sum_{(i,j)\in\A_k}Q_{i,j}^{(k)}(r)N^i
\Delta_\sigma^{\frac{j}{2}}u.$$
Then, with hypothesis $(b)$, for $0\leq j\leq\frac{k}{2}$,
$\mm_\alpha\ent{\Delta_\sigma^j u}\in L^p$ thus

$$\mm_\alpha\ent{\sum_{j=1}^{\ent{\frac{k}{2}}}P_j^{(k)}(r)\Delta_\sigma^j
u}
\in L^p(\S^{n-1}).$$

On the other hand, corollary 9 implies that

$$\mm_\alpha\ent{(1-r^2)\sum_{(i,j)\in\A_k}Q_{i,j}^{(k)}(r)N^i
\Delta_\sigma^{\frac{j}{2}}u}\leq
C\sum_{(i,j)\in\A_{k-1}\cup{(0,0)}}\mm_\beta\ent{N^i
\Delta_\sigma^{\frac{j}{2}}u}
\in L^p(\S^{n-1}),$$
since by proposition 16 $\norm{\mm_\alpha\ent{N^i\ffi}}_{L^p}$
and
$\norm{\mm_\beta\ent{N^i\ffi}}_{L^p}$ are equivalent for $\hh$-harmonic
functions and since, by induction hypothesis, $\mm_\alpha\ent{N^i
\Delta_\sigma^{\frac{j}{2}}u}\in L^p(\S^{n-1})$. We deduce from it
that $\mm_\alpha\ent{v}\in L^p(\S^{n-1})$. But, solving the differential
equation (5.2), we get

$$N^ku(r\zeta)=\left(\frac{1-r^2}{r^2}\right)^{n-1-k}
\int_0^rv(t\zeta)\frac{t^{2n-3-2k}}{(1-t^2)^{n-k}}dt
\leqno(5.3)$$
thus $\mm_\alpha\ent{N^ku}\in L^p(\S^{n-1})$ since $k\leq n-2$.

Assume now that $u$ satisfies $(a)$ {\it i.e.} that 
$\mm_\alpha\ent{N^ju}\in L^p(\S^{n-1})$ for $j\leq k$
and let us show that $\mm_\alpha\ent{\Delta_\sigma^ju}\in L^p(\S^{n-1})$
for $j\leq\frac{k}{2}$. Let $1>\beta>\alpha$.

According to the induction hypothesis, for $j\leq\frac{k}{2}-1$,
$\mm_\alpha\ent{\Delta_\sigma^ju}\in L^p(\S^{n-1})$ then

$$\mm_\alpha\ent{\sum_{j=1}^{\frac{k}{2}-1}P_j^{(k)}(r)
\Delta_\sigma^ju}\in L^p(\S^{n-1}).$$
With corollary 9, $\mm_\alpha\ent{(1-r^2)N^{k+1}u}\leq
C\mm_\beta\ent{N^ku}\in L^p(\S^{n-1})$, by proposition 16. 
As, for a regular function $\ffi$,

$$\mm_\alpha\ent{(1-r^2)\ffi}\leq C\bigl(\mm_\alpha\ent{(1-r^2)^2N\ffi}
+\abs{\ffi(0)}\bigr).$$
Finally, iterating this inequality, for $(i,j)\in\A_k$,

\begin{align}
\mm_\alpha\ent{(1-r^2)Q_{i,j}^{(k)}(r)N^i\Delta_\sigma^{\frac{j}{2}}u}
\leq&C\mm_\alpha\ent{(1-r^2)N^i\Delta_\sigma^{\frac{j}{2}}u}\notag\\
\leq&C\mm_\alpha\ent{(1-r^2)^{k-i+1}N^k\Delta_\sigma^{\frac{j}{2}}u}
+C\sum_{j=0}^{k-1}\abs{\nabla^ju(0)}\notag\\
\leq&C\sum_{j=0}^k\mm_\beta\ent{N^ju},\notag
\end{align}
by corollary 9 and $k-i+1\geq j$ since
$(i,j)\in\A_k$. So, with proposition 16,
$C\sum_{j=0}^k\mm_\beta\ent{N^ju}\in L^p(\S^{n-1})$, thus

$$\mm_\alpha\ent{\sum_{(i,j)\in\A_k}
(1-r^2)Q_{i,j}^{(k)}(r)N^i\Delta_\sigma^{\frac{j}{2}}u}\in
L^p(\S^{n-1}).$$
We then get from lemma 22 that
$\mm_\alpha\ent{P_{k/2}^{(k)}(r)\Delta_\sigma^{k/2}u}\in L^p(\S^{n-1})$,
and as $P_{k/2}^{(k)}$ is not zero on the boundary, 
$\mm_\alpha\ent{\Delta_\sigma^{k/2}u}\in L^p(\S^{n-1})$. So $(a)$ and
$(b)$ are equivalent.

Let us now show that $(a)+(b)$ implies $(c)$. It is enough to show this
implication for $\X$ a differential operator of the form $\X=N^j\Y$ with
$\Y$
a product of $k-j$ operators of the form $\lll_{i,j}$. 
We can assume that $j<k$. Let

$$v=(1-r^2)N^{j+1}u+2(n-1-j)N^ju$$
and compose with $\Y$, it results that

$$\Y v=(1-r^2)N^{j+1}\Y u+2(n-1-j)N^j\Y u.$$
Using as previously formula (5.3), we see that

$$\mm_\alpha\ent{\X u}\in L^p(\S^{n-1})$$
which completes the proof in the case $k$ is even.

Assume now that $k$ is {\it odd}. The proof of $(b)\Rightarrow(a)$ is
similar to case $k$ even. The converse is again based on lemma 22.
According to the induction hypothesis,
$\mm_\alpha\ent{\Delta_\sigma^lu}\in L^p(\S^{n-1})$ for $0\leq
l\leq\frac{k-1}{2}=\ent{\frac{k}{2}}$ so that

$$\mm_\alpha\ent{\sum_{j=1}^{\ent{\frac{k}{2}}}P_j^{(k)}(r)\Delta_\sigma^ju}\in
L^p(\S^{n-1}).$$
One has, as before,

$$\mm_\alpha\ent{(1-r^2)\sum_{(i,j)\in\A_k\setminus(0,k+1)}Q_{i,j}^{(k)}(r)
N^i\Delta_\sigma^{\frac{j}{2}}u}\in L^p(\S^{n-1})$$
and that $\mm_\alpha\ent{(1-r^2)N^{k+1}}\in L^p(\S^{n-1})$.

Combining all this, we get that

$$\mm_\alpha\ent{(1-r^2)Q_{0,k+1}^{(k)}(r)\Delta_\sigma^{\frac{k+1}{2}}u}\in
L^p(\S^{n-1})$$
and as $Q_{0,k+1}^{(k)}$ is non-zero on the boundary, we finaly get
$$\mm_\alpha\ent{(1-r^2)\Delta_\sigma^{\frac{k+1}{2}}u}\in
L^p(\S^{n-1})$$
and $(a)$ and $(b)$ are equivalent.\hfill$\Box$
\end{proof}

We will now prove the area integral characterization
in theorem B.

\begin{proof}[Proof of theorem B] The fact that $S_\alpha$ can be
replaced by $S_\alpha^N$ is a direct consequence of
lemma 14\ (with $\gamma=-\frac{n}{2}+1$). 
Further, as

$$\abs{N\Delta_\sigma^{k/2}u}\leq CI_k(N^{k+1}\Delta_\sigma^{k/2}u)
+\sup_{0\leq j\leq 2k,\abs{z}\leq\eps}\abs{\nabla^ju},$$
so, lemma 14 and the mean value inequality imply that

$$S_\alpha^N\ent{\Delta_\sigma^{k/2}u}\leq
S_\beta^N\ent{N^ku}+\norm{u}_{\hh^p},$$
so that $8$ implies $5$.

Let us now prove that if, for $0\leq j\leq\frac{k}{2}$,
$S_\alpha\ent{\Delta_\sigma^ju}\in L^p(\S^{n-1})$, then
$S^N_\alpha\ent{N^ku}\in L^p(\S^{n-1})$. The proof goes according to the
method developped for the equivalence of maximal functions.

For simplicity, we will restrict our attention to the case $k=1$. In
order to estimate $S_\alpha^N\ent{Nu}$, we have to estimate $Nu$.
While trying to use the previous method, lemma 22 for $k=2$ 
does not give a satisfying estimate. However, we can obtain the
desired estimate as follows. Denote by $v$ the function

$$v=2(n-2)Nu+(1-r^2)N^2u,$$
then $Nv=2(n-3)N^2u+(1-r^2)N^3u+2(1-r^2)N^2u$ and write this in the form
$Nv=w+2(1-r^2)N^2u$. As before, solving the differential equation  
$(1-r^2)N^3u+2(n-3)N^2u=w$, we have

$$N^2u(r\zeta)=\left(\frac{1-r^2}{r^2}\right)^{n-3}
\int_0^rw(tz)t^{2(n-3)+1}(1-t^2)^{2-n}dt$$
so that

$$\abs{N^2u(r\zeta)}\leq C(1-r^2)^{n-3}I_{3-n}(\abs{w}).
\leqno(5.4)$$
On the other hand, by lemma 22 for $k=1$,

$$v=(n-2)(1-r^2)Nu-(1-r^2)\Delta_\sigma u,$$
so

$$Nv=(n-2)(1-r^2)N^2u-(1-r^2)N\Delta_\sigma u+2r^2\Delta_\sigma
u-2r^2(n-2)Nu.$$

Recall that $\abs{f}\leq CI_l(\abs{N^lf})+C\sup_{\abs{z}<\eps,j\leq l}
\abs{\nabla^jf}$, and that $I_{k+1}(\abs{f})\leq I_k(\abs{f})$. Using these
facts, one gets

\begin{align}
\abs{w}\leq&\abs{Nv}+2(1-r^2)\abs{N^2u}\notag\\
\leq&C\left(I_1\bigl(\abs{N\Delta_\sigma u}\bigr)
+I_1\bigl(\abs{N^2u}\bigr)+(1-r^2)\abs{N^2u}
+(1-r^2)\abs{N\Delta_\sigma u}\right)+\sup_{0\leq j\leq
3,\abs{z}\leq\eps}\abs{\nabla^ju}.\notag
\end{align}
Inserting this in (5.4), and invoking the facts that
$I_l\bigl((1-r^2)\abs{f}\bigr)=I_{l+1}(\abs{f})$ and that
$I_l(I_s\abs{f})\leq CI_{l+s}(\abs{f})$, one gets

$$\abs{N^2u(r\zeta)}\leq
C(1-r^2)^{n-3}\left(I_{4-n}(\abs{N\Delta_\sigma
u})+I_{4-n}(\abs{N^2u})+\sup_{\abs{z}<\eps,0\leq j\leq
3}\abs{\nabla^ju}\right).$$

We are now in position to estimate $S_\alpha^N\ent{Nu}$ :

\begin{align}
S_\alpha^N\ent{Nu}(\zeta)^2=&\int_{\aa_\alpha(\zeta)}
\abs{N^2u(x)}^2(1-\abs{x})^{2-n}dx\notag\\
\leq&C\left(\int_{\aa_\alpha(\zeta)}\ent{I_{4-n}(\abs{N\Delta_\sigma
u})}^2(1-\abs{x}^2)^{n-4}dx\right.\notag\\
&\left.+\int_{\aa_\alpha(\zeta)}\ent{I_{4-n}(\abs{N^2u})}^2
(1-\abs{x}^2)^{n-4}dx+\sup_{\abs{z}<\eps,0\leq j\leq
3}\abs{\nabla^ju}^2\right)\notag
\end{align}
A further appeal to lemma 13, with $l=4-n$, $d=0$, $k=2$ and
$\gamma=\frac{n-4}{2}$ leads to

$$\int_{\aa_\alpha(\zeta)}\ent{I_{4-n}(\abs{N^2u})}^2
(1-\abs{x}^2)^{n-4}dx\leq C\int_{\aa_\beta(\zeta)}
\abs{N^2 u}^2(1-\abs{x}^2)^{4-n}dx.$$

A last appeal to lemma 13, with $l=4-n$, $d=0$, $k=1$
and $\gamma=\frac{n-4}{2}$ leads to

$$\int_{\aa_\alpha(\zeta)}\ent{I_{4-n}(\abs{N\Delta_\sigma
u})}^2(1-\abs{x}^2)^{n-4}dx\leq C\int_{\aa_\beta(\zeta)}
\abs{N\Delta_\sigma u}^2(1-\abs{x}^2)^{4-n}dx
\leq CS_\gamma^N\ent{\Delta^{\frac{1}{2}}_\sigma u},$$
by the mean value properties. As the only part that matters in this last integral is the part near to
the boundary, we will cut it into two parts. Let $\kappa$ be a constant that
we will fix later. Then

\begin{align}\int_{\aa_\beta(\zeta)}
\abs{N^2 u}^2(1-\abs{x}^2)^{4-n}dx\leq&
\int_{\aa_\beta(\zeta)\cap B(0,\kappa)}\abs{N^2u}^2
(1-\abs{x}^2)^{4-n}dx\notag\\
&+\int_{\aa_\beta(\zeta)\cap\bigl(\B_n\setminus B(0,\kappa)\bigr)}
\abs{N^2u}^2(1-\abs{x}^2)^{4-n}dx\notag\\
\leq&C\sup_{\abs{z}\leq\kappa,0\leq j\leq3}\abs{\nabla^ju}
+(1-\kappa^2)^2\int_{\aa_\beta(\zeta)}
\abs{N^2 u}^2(1-\abs{x}^2)^{2-n}dx.\notag
\end{align}
Grouping the above estimates, we finaly get

$$S_\alpha^N\ent{Nu}(\zeta)^2\leq
CS_\gamma^N\ent{\Delta_\sigma^{\frac{1}{2}}u}(\zeta)^2
+C(1-\kappa^2)^2S_\beta\ent{Nu}(\zeta)^2
+C\sup_{\abs{z}\leq\kappa,0\leq
j\leq3}\abs{\nabla^ju}.\leqno(5.5)$$

But this inequality depends only on the mean value inequality, in
particular, one can replace $u$ in (5.5) by 
$u_\delta(x)=u(\delta x)$ and get

$$S_\alpha^N\ent{Nu_\delta}(\zeta)^2\leq 
CS_\gamma^N\ent{\Delta_\sigma^{\frac{1}{2}}u_\delta}(\zeta)^2
+C(1-\kappa^2)^2S_\beta^N\ent{Nu_\delta}(\zeta)^2
+C\sup_{\abs{z}\leq\kappa,0\leq j\leq3}\abs{\nabla^ju_\delta}$$
with constants independant on $\frac{1}{2}<\delta<1$. Then taking
$L^p(\S^{n-1})$ norms, one gets

\begin{align}
\norm{S_\alpha^N\ent{Nu_\delta}}_p\leq&
C\norm{S_\gamma^N\ent{\Delta_\sigma^{\frac{1}{2}}u_\delta}}_p
+C(1-\kappa^2)\norm{S_\beta^N\ent{Nu_\delta}}_p+C\norm{u}_{\hh^p}\notag\\
\leq&C\norm{S_\gamma^N\ent{\Delta_\sigma^{\frac{1}{2}}u}}_p
+C'(1-\kappa^2)\norm{S_\alpha^N\ent{Nu_\delta}}_p+C\norm{u}_{\hh^p},\notag
\end{align}
with lemma 12. It is now enough to choose $\kappa$ such that
$C'(1-\kappa^2)=\frac{1}{2}$, then

$$\norm{S_\alpha^N\ent{Nu_\delta}}_p\leq
\frac{C}{2}\norm{S_\gamma^N\ent{\Delta_\sigma^{\frac{1}{2}}u}}_p+
\frac{C}{2}\norm{u}_{\hh^p}.$$
We then conclude by monotone convergence when $\delta\rightarrow1$.

So far we have proved the equivalence of properties $1$ to $8$ and that
$9$ implies these properties. To see that $9$ is actually equivalent to
them, it is enough to see that in $9$, $\nabla^k$ can be replaced by any
operator of the form $N\Y$ where $\Y$ is a product of $k-1$ operators
of the form $\lll_{i,j}$. This is then a direct consequence of $2$ and
$7$.\hfill$\Box$
\end{proof}

\section{Lipschitz spaces and Zygmund classes}

\subsection{Lipschitz spaces}\ \\

\begin{lemmanum}[24]Let $k\in\N$ and $0<\alpha<1$. Assume further that
if $n$ is odd then $k\leq n-2$. There exists a constant $C$ such that
for every $f\in\cc^{k+\alpha}(\S^{n-1})$, for every $r\xi\in\B_n$,

$$(1-r^2)^k\abs{\nabla^k\P_h\ent{f}(r\xi)}\leq C(1-r^2)^\alpha.$$
In particular, $\P_h\ent{f}\in\cc^{k+\alpha}(\overline{\B_n})$.
\end{lemmanum}

\begin{proof}[Proof] Fix $\xi_0\in\S^{n-1}$, there exists
$P_{\xi_0}^{(k)}$, a combination of
spherical harmonics of order less than $k$  such that the
Taylor polynomials of order $k$ at $\xi_0$ of $P_{\xi_0}^{(k)}$ and of
$f$ coincide. Then

$$\abs{f(\xi)-P_{\xi_0}^{(k)}(\xi)}\leq C(1-<\xi,\xi_0>)^\alpha.$$
But then,
$\P_h\ent{f}=\P_h\ent{P_{\xi_0}^{(k)}}+\P_h\ent{f-P_{\xi_0}^{(k)}}$.

From the spherical harmonics expansion of $\P_h$, if $n$ is {\it odd}
and $k\leq n-2$, or if $n$ is {\it even}, there exists a constant $C$,
independent of $\xi_0$ such that

$$\norm{\P_h\ent{P_{\xi_0}^{(k)}}}_{\cc^k}\leq C\norm{f}_{\cc^k}.$$
To estimate $\nabla^k\P_h\ent{f-P_{\xi_0}^{(k)}}$ we need the following
estimates on the hyperbolic Poisson kernel :

\begin{enumerate}
\item $\displaystyle\abs{\nabla^k\P_h(r\zeta,\xi)}\leq\frac{C}{(1-r)^{n-1+k}},$

\item
$\displaystyle\abs{\nabla^k\P_h(r\zeta,\xi)}\leq
\frac{C}{(1-<\zeta,\xi>)^{n-1+k}}$ provided
$r\zeta\in\aa_{\alpha_k}(\xi)$ for some $\alpha_k$ small enough.
\end{enumerate}

Both estimates result directly from the mean value inequalities applied
to the $\hh$-harmonic function $u(r\zeta)=\P_h(r\zeta,\xi)$,
$\xi\in\S^{n-1}$ fixed.

Furthermore, we are only interested in the estimates when $r\zeta$ is
``near'' to $\S^{n-1}$, {\it i.e.} $r$ near to $1$. In this case $(2)$ holds
when $1-\scal{\zeta,\xi}<c_k(1-r)$. Then

\begin{align}
(1-r)^k\abs{\nabla^k\P_h\ent{f-P_{\xi_0}^{(k)}}(r\xi_0)}\leq&
(1-r)^k\int_{\S^{n-1}}\abs{\nabla^k\P_h(r\zeta,\xi_0)}
\abs{f(\xi)-P_{\xi_0}^{(k)}}d\sigma(\xi)\notag\\
\leq&(1-r)^k\int_{1-\scal{\xi,\xi_0}<c_k(1-r)}\abs{\nabla^k\P_h(r\zeta,\xi_0)}
\abs{f(\xi)-P_{\xi_0}^{(k)}}d\sigma(\xi)\notag\\
&+(1-r)^k\int_{1-\scal{\xi,\xi_0}>c_k(1-r)}\abs{\nabla^k\P_h(r\zeta,\xi_0)}
\abs{f(\xi)-P_{\xi_0}^{(k)}}d\sigma(\xi)\notag\\
\leq&\frac{C}{(1-r)^{n-1}}\int_{1-\scal{\xi,\xi_0}<c_k(1-r)}
(1-<\xi,\xi_0>)^\alpha d\sigma(\xi)\notag\\
&+C\int_{1-\scal{\xi,\xi_0}>c_k(1-r)}
(1-<\xi,\xi_0>)^{-(n-1)+\alpha}d\sigma(\xi)\notag
\end{align}
where for the first integral we have used estimate $(1)$ on $\P_h$ and
for the second we have used estimate $(2)$. This immediatly leads to the
desired result.
\hfill$\Box$
\end{proof}

\begin{remark} When $n$ is even, lemma 3 and the result in the
Euclidean case give directly the result.
\end{remark}

The converse of this result can also be obtained by a transfer from 
the Euclidean case :

\begin{lemmanum}[25]Let $k\in\N$ and $0<\alpha<1$. Let $f$ be a
distribution on $\S^{n-1}$ and let $u=\P_h\ent{f}$. Assume that
there exists a constant $C$ such that $u$ satisfies the
following inequality :

$$(1-r^2)^k\abs{\nabla^ku(r\xi)}\leq C(1-r^2)^\alpha.$$
Then $f\in\cc^{k+\alpha}(\S^{n-1})$.
\end{lemmanum}

\begin{proof}[Proof] Let $v=\P_e\ent{f}$. Then, by lemma 4, for
$r\zeta\in\B_n$,

$$v(r\zeta)=\int_0^1\eta(r,s)u(rs\zeta)ds.$$
But the estimates on $\eta$ imply that

$$(1-r^2)^k\abs{\nabla^kv(r\xi)}\leq C(1-r^2)^\alpha.$$
Thus, with the result on Euclidean harmonic functions,
$f\in\cc^{k+\alpha}(\S^{n-1})$.\hfill$\Box$
\end{proof}

\begin{remark} When $n$ is odd, if $k\geq n-1$, the condition
$(1-r^2)^k\abs{\nabla^ku(r\xi)}\leq C(1-r^2)^\alpha.$ is reduced to $u$
constant.
\end{remark}

\subsection{Zygmund classes}
Let us fix $\xi_0\in\S^{n-1}$ and denote
for $\xi\in\S^{n-1}$ by $R_\xi$ a rotation on $\S^{n-1}$ that maps
$\xi_0$ to $\xi$. Let $R_\xi^*$ be the reverse rotation that maps $\xi$
to $\xi_0$.

Define the Zygmund class of order $n$ on $\S^{n-1}$ by

$$\Z_n(\S^{n-1})=\{f\in\cc^{n-2}(\S^{n-1})\ :\
\abs{\tilde\Delta_\xi^nf(\zeta)}\leq C(1-\scal{\zeta,\xi})^{n-1}\}$$
where $\tilde\Delta_\xi^n$ is the difference operator defined by induction
on $j$ by

$$\tilde\Delta_\xi^1f(\zeta)=f(R_\xi\zeta)-f(\zeta)$$
and $\tilde\Delta_\xi^{j+1}f=\tilde\Delta_\xi^1(\tilde\Delta_\xi^jf)$.

Define the Zygmund class of order $n$ on $\B_n$ by

$$Z_n(\B_n)=\{f\in\cc^{n-2}(\B_n)\ :\
\norm{\Delta_h^{n,\gamma}f(\zeta)}_{L^\infty(\B_n)}\leq C\abs{h}^{n-1}
\ \mathrm{for\ any\ curve\ }\gamma:\ent{0,1}\mapsto\B_n\}$$
where $\Delta_h^j$ is the difference operator along $\gamma$ also
defined inductively by

$$\Delta_h^{1,\gamma}u(z)=u\circ\gamma(h)-u\circ\gamma(0)$$
where $\gamma(0)=z$, and $\Delta_h^{j+1,\gamma}u=\Delta_h^{1,\gamma}
(\Delta_h^{j,\gamma}u)$.

It follows from standard methods (using the mean-value properties) that
the set of $\hh$-harmonic functions belonging to the Zygmund class
is given by :

\begin{align}
\{u\ \hh-\mathrm{harmonic},\ u\in Z_n(\B_n)\}
=&\left\{u\in\cc^{n-1}(\B_n),\ u\ \hh-\mathrm{harmonic},\
\abs{\nabla^nu(z)}\leq\frac{C}{1-\abs{z}}\right\}\notag\\
=&\Bigl\{u\in\cc^{n-1}(\B_n),\ u\ \hh-\mathrm{harmonic},\ 
\abs{N^nu(z)}\leq\frac{C}{1-\abs{z}}\Bigr\}.\notag
\end{align}
The next theorem states that this class is the set of $\hh$-harmonic
extensions of members of the Zygmund class of order $n$ on $\S^{n-1}$.

\begin{theorem}[26] A function $f$ belongs to $\Z_n(\S^{n-1})$ if
and only if $u=\P_h\ent{f}$ belongs to $\Z_n(\B_n)$.
\end{theorem}

\begin{proof}[Proof] The proof follows from standard arguments. Let us
first prove that $\P_h\ent{f}\in\Z_n(\B_n)$ when $f\in\Z_n(\S^{n-1})$.
Note that, for fixed $\xi\in\S^{n-1}$,

$$\int_{\S^{n-1}}N^j_z\P_h(z\zeta,\xi)d\sigma(\zeta)=0$$
for any $j\geq 1$, since $\P_h$ has integral $1$ on $\S^{n-1}$. Using
this fact and the symmetry under rotations of $\P_h(r\zeta,.)$, we
get that

$$\int_{S^{n-1}}\P_h(r\zeta,\xi)f(\xi)d\sigma(\xi)=
\frac{1}{n}\int_{\S^{n-1}}N^n\P_h(r\zeta,\xi)\tilde\Delta_\xi^n
f\bigl((R_\xi^*)^{n-1}\xi_0\bigr)d\sigma(\xi).$$
Now, by assumption, $\abs{\tilde\Delta_\xi^n
f\bigl((R_\xi^*)^{n-1}\xi_0\bigr)}\leq C(1-\scal{\xi,\xi_0})^{n-1}$, so
that the desired estimate follows as in the proof of lemma 24.

For the converse, by the same proof as in lemma 25, we get that
$v=\P_e\ent{f}$ belongs to $\Z_n(\B_n)$ and we conclude that
$f\in\Z_n(\S^{n-1})$ from the euclidean harmonic theory.
\end{proof}

\begin{remark} It is proved in \cite{Ja3} that any $\hh$-harmonic
function $u$ is at most in $\Z_n(\B_n)$. In other words, it means that
the Zygmund class of order $n$ is the limit class preserved by the
hyperbolic Poisson kernel.
\end{remark}

\bibliographystyle{plain}
\bibliography{harm}
\end{document}